\documentclass[10pt]{article}
\usepackage{mathrsfs}
\usepackage{amsfonts}
\usepackage{amsthm,amsmath,amssymb,anysize}
\newtheorem{lemma}{Lemma}[section]
\newtheorem{theorem}[lemma]{Theorem}
\newtheorem{remark}[lemma]{Remark}
\newtheorem{proposition}[lemma]{Proposition}

\newtheorem{corollary}[lemma]{Corollary}
\usepackage[numbers,sort&compress]{natbib}
\marginsize{3.6cm}{3.6cm}{3.6cm}{3cm}
\numberwithin{equation}{section}

\title{\textsf{Finite dimensional  special odd contact
superalgebras over a field of prime characteristic}}
\author{\textsc{Wende Liu$^{1, 2,}$}\footnote{Supported by the NSF
  of China (10871057) and the NSF of Heilongjiang Province (A200802)}\;\;\textsc{and Jixia Yuan$^{1,}$}\footnote{Correspondence:
  jxy@hrbnu.edu.cn (J. Yuan); wendeliu@ustc.edu.cn (W. Liu)}   \\
  \\
  \textit{$^{1}$Department of Mathematics},
  \textit{Harbin Institute of Technology}\\
  \textit{Harbin 150006, China}\\
\\
  \ \ \textit{$^{2}$School of Mathematical Sciences},
  \textit{Harbin Normal University} \\
  \textit{Harbin 150025, China}
  }
\date{ }
\begin{document}
\maketitle
\begin{quotation}
\noindent\textbf{Abstract} This paper considers  a family of
finite dimensional simple  Lie superalgebras of Cartan type over a
field of characteristic $p>3$, the so-called special odd contact
superalgebras. First, the spanning sets are determined for the Lie superalgebras and their relatives.
Second,  the spanning sets are used to characterize the
simplicity and to compute the dimension formulas. Third, we determine the
superderivation algebras and the first cohomology groups.
   Finally,  the dimension formulas
and the first cohomology groups are
used to make a comparison between the special odd contact
superalgebras and the other simple  Lie superalgebras of Cartan
type.\\



  \end{quotation}

  \setcounter{section}{-1}
\section{Introduction}

\noindent The theory of Lie superalgebras has seen a significant
development. For example, V. G. Kac classified the finite
dimensional simple Lie superalgebras and the infinite dimensional
simple linearly compact Lie superalgebras over algebraically closed
fields of characteristic zero (see \cite{k1,k2,sc}). But there are
not so plentiful results for modular Lie superalgebras (that is, Lie
superalgebras over fields of prime characteristic). The
classification problem is still
 open for finite dimensional simple modular Lie superalgebras (see \cite{bl,z} for example). As far as we know,
 \cite{kl,p} should be two of the earliest
 papers on modular Lie superalgebras.  Recently, certain new simple
 Lie superalgebras over a field of  characteristic 3
 were constructed and studied \cite{bl,e1}.
 In
 \cite{fzj,lzw,z}
    six families of finite dimensional
  modular Lie superalgebras of Cartan type $W,S,H,K,HO$ and $KO$
 were considered and the simplicity and
 restrictiveness were determined.
The superderivation  algebras have been sufficiently studied
for these Lie superalgebras \cite{lz3,lzw,wz,zz}. In a recent paper \cite{lh}, the finite dimensional special
 odd Hamiltonian superalgebras were introduced and the spanning set,
  simplicity and dimension formula were determined.

 In the present paper, motivated by  \cite{k2,lh}, we study a  family
of finite dimensional  Lie superalgebras of Cartan type over a filed
of characteristic $p>3,$ the special odd contact superalgebras, and
determine the spanning sets, simplicity, dimension formulas, the
superderivations and the first cohomology.   From the discussions and certain conclusions in this paper one may  see that, as in the non-modular case,
 the special odd contact superalgebras possess more complicated structures and
have no analogous in the finite dimensional modular Lie algebras of
Cartan type (cf. \cite{st,sf}).

\section{Preliminaries}

Throughout $\mathbb{F}$ is a field of characteristic $p>3$;
$\mathbb{Z}_2:= \{\bar{0},\bar{1}\}$ is the additive group of two
elements; $\mathbb{N}$ and $\mathbb{N}_0$ are the sets of positive
integers and nonnegative integers, respectively. Fix an integer
$n\geq 3$ and an $n$-tuple $\underline t:=(t_1,\ldots,t_n)\in
\mathbb{N}^n$. Put $\pi:=(\pi_1,\ldots,\pi_n),$ where
$\pi_i:=p^{t_i}-1$ for $i\in \overline{1,n}.$ Let
$\mathcal{O}(n;\underline{t})$ be the divided power algebra  with
$\mathbb{F}$-basis $\{x^{(\alpha)}\mid \alpha\in
\mathbb{A}(n;\underline{t})\}$, where
$\mathbb{A}(n;\underline{t}):=\{\alpha\in \mathbb{N}_0^n\mid
\alpha_i\leq \pi_i\}$. Note that
$x^{(0)}:=1\in\mathcal{O}(n;\underline{t})$, where
$0=(0,\ldots,0)\in \mathbb{A}(n;\underline{t}).$ For
$\varepsilon_i:=(\delta_{i1},\delta_{i2},\ldots,\delta_{in})\in
\mathbb{A}(n;\underline{t})$,
  write $x_i$ for $x^{(\varepsilon_i)}$, where $i=\overline{1,n}.$ Let
$\Lambda(m)$ be the exterior superalgebra over $\mathbb{F}$ in $m$
variables $x_{n+1},x_{n+2},\ldots,x_{n+m}$.
 Set
$$\mathbb{B}(m):=\left\{ \langle i_1,i_2,\ldots,i_k\rangle \mid n+1\leq i_1<i_2<\cdots <i_k\leq n+m;\;
k\in \overline{0, m}\right\}.$$
 For
 $u:=\langle i_1,i_2,\ldots,i_k\rangle \in \mathbb{B}(m),$ write $ |u| :=k$
   $x^u:=x_{i_1}x_{i_2}\cdots x_{i_k},$ and
  denote by $u$ itself the index set $\{i_1,i_2,\ldots,i_k\}$.
  For $u,\upsilon \in \mathbb{B}(m)$ with $u\cap
  \upsilon=\emptyset,$ write $u+\upsilon$ for the uniquely determined element
  $w\in \mathbb{B}(m)$ such that $ w=u\cup
  \upsilon.$ If $\upsilon\subset u,$ write $u-\upsilon$ for the uniquely determined element $w\in
  \mathbb{B}(m)$ such that $w=u\setminus\upsilon.$
  Clearly, the associative superalgebra
$$\mathcal{O}(n,m;\underline{t}):=\mathcal{O}(n;\underline{t})\otimes
\Lambda(m)$$ has a so-called \textit{standard}
$\mathbb{F}$-\textit{basis} $$\{ x^{(\alpha ) }\otimes x^u\mid
(\alpha,u) \in \mathbb{A}(n;\underline{t})\times \mathbb{B}(m)\}.$$

 Let
$\partial_r$ be the superderivation of
$\mathcal{O}(n,m;\underline{t}) $ defined by
$\partial_{r}(x^{(\alpha)})=x^{(\alpha-\varepsilon_{r})}$ for $r\in
\overline{1,n} $ and $\partial_{r}(x_{s})=\delta_{rs}$ for $r,s\in
\overline{1, n+m}.$ The \textit{generalized Witt superalgebra} $
W\left(n,m;\underline{t}\right)$ is  $\mathbb{F}$-spanned by
$\{f_r
\partial_r\mid f_r\in \mathcal{O}(n,m;\underline{t}), r\in
\overline{1,n+m}\}.$   Note that $W(n,m;\underline{t}) $ is a free $
\mathcal{O} \left(n,m;\underline{t}\right)$-module with basis $ \{
\partial_r\mid r\in
 \overline{1,n+m}\}.$ In particular, $W (n,m;\underline{t} ) $ has a
\textit{standard} $\mathbb F$-\textit{basis}
 $$\{x^{(\alpha)}x^{u}\partial_r\mid (\alpha,u,r)\in
\mathbb{A}(n;\underline{t})\times\mathbb{B}(m)\times\overline{1,n+m}\}.$$

 For an $n$-tuple
$\alpha:=(\alpha_1,\ldots,\alpha_n)\in \mathbb{N}_0^n$, put
$|\alpha|:=\sum_{i=1}^n\alpha_i$. The
associative algebra $\mathcal{O}(n,m;\underline t)$  has a standard
$\mathbb{Z}$-grading structure $\mathcal{O}(n,m;\underline
t)=\oplus_{i=0}^{\xi}\mathcal{O}(n,m;\underline t)_{i}$, where
$$\mathcal{O}(n,m;\underline t)_{i}:={\rm
span}_{\mathbb{F}}\{x^{(\alpha)}x^u\mid |\alpha|+|u|=i\},\quad
\xi:=|\pi|+m.
$$
 For a   vector superspace $V=V_{\bar{0}}\oplus
V_{\bar{1}}, $ we write $\mathrm{p}(x):=\theta$ for the
\textit{parity} of a homogeneous element $x\in V_{\theta}, $
$\theta\in \mathbb{Z}_{2}.$ Once the symbol $\mathrm{p}(x)$ appears
in this paper, it
 will imply that $x$ is  a $\mathbb{Z}_2$-homogeneous element.

When $m=n+1$, we write $\mathcal {O}:=\mathcal
{O}(n,n+1;\underline{t}).$
 Recall the odd contact superalgebra, which is a simple Lie superalgebra contained in $W(n,n+1;\underline{t}),$
defined as follows (see \cite{fzj,k2}):
$$
KO(n,n+1;\underline{t}):=\{D_{KO}(a)\mid a\in \mathcal{O}\},
$$
where
$$D_{KO}(a):=T_{H}(a)+(-1)^{\mathrm{p}(a)}\partial_{2n+1}(a)\mathfrak{D}+ (\mathfrak{D}(a)-2a )\partial_{2n+1}$$
and $$\mathfrak{D}:=\sum _{i=1}^{2n}x_{i}\partial_{i},\quad
T_{H}(a):=\sum
_{i=1}^{2n}(-1)^{\mu(i')\mathrm{p}(a)}\partial_{i'}(a)\partial_{i},$$
\[i':= \left\{\begin{array}{ll}
i+n, &\mbox{if}\; i\in\overline{1, n}
\\i-n, &\mbox{if}\; i\in \overline{n+1, 2n},
\end{array}\right.\quad
\mu(i):=\left\{\begin{array}{ll}\bar{0},  &\mbox{if}\; i\in \overline{1,n}\\
\bar{1},&\mbox{if}\; i\in\overline{n+1,2n+1}.
\end{array}\right. \] Note that  for $a,b\in \mathcal{O} $ ,
\begin{equation}\label{liuee1}
[D_{KO}(a),
D_{KO}(b)]=D_{KO}(D_{KO}(a)(b)-(-1)^{\mathrm{p}(a)}2\partial_{2n+1}(a)b)
\end{equation}
(see \cite{fzj,k2}). Given $\lambda\in \mathbb{F},$ put
$$
SKO''(n,n+1;\lambda,\underline{t}):=\{D_{KO}(a)\mid
\mathrm{div}_{\lambda}(a)=0, a\in \mathcal{O}\},
$$
where
$$
\mathrm{div}_{\lambda}(a):=(-1)^{\mathrm{p}(a)}2(\Delta(a)+(\mathfrak{D}-n\lambda\mathrm{id}_{\mathcal{O}})
\partial_{2n+1}(a))
$$
 and
 $$\Delta:=\sum_{i=1}^{n}\partial_{i}\partial_{i'}.$$
Then one can verify that $SKO''(n,n+1;\lambda,\underline{t})$ is a
subalgebra of $KO(n,n+1;\lambda,\underline{t})$ (c.f. \cite{k2}).
The following symbols will be frequently used in this paper:
$$\Delta_{i}:=\partial_{i}\partial_{i'};
\nabla_{i}(x^{(\alpha)}x^{u}):=x^{(\alpha+\varepsilon_{i})}x_{i'}x^{u};
\Gamma_{i}^{j}:=\nabla_{j}\Delta_{i}\quad \mbox{for}\; i,j\in
\overline{1,n},$$  for fixed $(\alpha,u) \in
\mathbb{A}(n;\underline{t})\times \mathbb{B}(n+1),$
$$
I(x^{(\alpha)}x^{u}):=I(\alpha,u):=\{i\in  \overline{1, n}\mid
\Delta_{i}(x^{(\alpha)}x^{u})\neq0\};
$$
$$
\widetilde{I}(x^{(\alpha)}x^{u}):=\widetilde{I}(\alpha,u):=\{i\in
\overline{1,n}\mid \nabla_{i}(x^{(\alpha)}x^{u})\neq0\}
$$
and
$$
\mathcal{D}^{\ast}:=\{x^{(\alpha)}x^{u}\in \mathcal {O}\mid
I(\alpha,u)\neq\emptyset,\widetilde{I}(\alpha,u)\neq\emptyset\}.
$$
Recall some Lie superalgebras contained in
$W(n,n;\underline{t})$:
\begin{eqnarray}\label{e1.1}
&& SHO'(n,n;\underline{t}):=\{T_{H}(a)\mid a\in \mathcal
{O}(n,n;\underline{t}), \Delta(a)=0\};\nonumber
\\
&& \overline{SHO}(n,n;t):=[SHO'(n,n;t),SHO'(n,n;t)].\nonumber
\end{eqnarray}
It was proven in \cite{lh} that
\begin{equation}
 \overline{SHO}(n,n;t) =\mathrm{span}_{\mathbb{F}}\big(\big\{
T_{H}\big(x^{(\alpha)}x^{u}\big)\mid I(\alpha,u)=\emptyset,
\widetilde{I}(\alpha,u)\neq \emptyset\big\}\cup\mathcal {G}\big),
\end{equation}
where
$$
\mathcal {G}:=\Big\{ T_{H}\Big(x^{(\alpha)}x^{u}-\sum _{i\in
I(\alpha,u)}\Gamma_{i}^{q}(x^{(\alpha)}x^{u})\Big)\mid
x^{(\alpha)}x^{u}\in \mathcal{D}^{\ast}, q\in
\widetilde{I}(\alpha,u)\Big\}.
$$
Let $\mathcal{K}(n,n+1;\lambda,\underline{t})$ be the subspace
spanned by the elements $D_{KO}(a),$ where $a\in \oplus_{i\geq
1}\mathcal {O}_{i}$ such that
\begin{equation}\label{liue1}
  \mathrm{div}_{\lambda}(a)=0,\;\partial_{2n+1}(a)=0.
\end{equation}
If $\partial_{2n+1}(a)=\partial_{2n+1}(b)=0,$ that is, $a, b\in
\mathcal {O}(n,n;\underline{t}),$ then by (\ref{liuee1}), we have
\begin{equation*}\label{hee1.1}
   [D_{KO}(a),D_{KO}(b)]=D_{KO}(T_{H}(a)(b)).
\end{equation*}
Note that
\begin{equation*}\label{hee1.1}
   [T_{H}(a),T_{H}(b)]=T_{H}(T_{H}(a)(b)).
\end{equation*}
It follows that $\mathcal{K}(n,n+1;\lambda,\underline{t})$ is a
subalgebra of $SKO'(n,n+1;\lambda,\underline{t}).$ Moreover, the
mapping
\begin{equation}\label{liue2}
\rho: \mathcal{K}(n,n+1;\lambda,\underline{t})\longrightarrow
SHO'(n,n;\underline{t}),\quad D_{KO}(a)\longmapsto T_{H}(a)
\end{equation}
 is an isomorphism of Lie superalgebras.
In this paper we mainly study  the following derived superalgebras
$$SKO'(n,n+1;\lambda,\underline{t}):=[SKO''(n,n+1;\lambda,\underline{t}), SKO''(n,n+1;\lambda,\underline{t})],$$
$$SKO(n,n+1;\lambda,\underline{t}):=[SKO'(n,n+1;\lambda,\underline{t}), SKO'(n,n+1;\lambda,\underline{t})],$$
called the \textit{special odd contact (Lie) superalgebras}. For
simplicity, in the following sections we shall write $\mathfrak{g}$
for $SKO$ and usually omit the parameters
$(n,n+1;\lambda,\underline{t}).$

\section{Spanning sets}
From now on,  we take the convention that the expression
$x^{(\alpha)}x^{u}$ implies that
 $\alpha\in \mathbb{A}(n;\underline{t})$ and $ u\in \mathbb{B}(n).$  For $f\in \mathcal {O},$ if there is $i\in \overline{1,n}$ such that
$\nabla_{i}(f)\neq 0$, then $f$ is called $i$-\textit{integral}.
 If $\widetilde{I}(\alpha,u)\neq \emptyset,$
we write $q_{(\alpha,u)}:=\mathrm{min}\widetilde{I}(\alpha,u).$
Recall our convention that $n\geq3$.  To formulate
 the linear generators of $\frak{g}''$ (Theorem \ref{t5}), we introduce the following symbols for $q\in \widetilde{I}(\alpha,u),$
\begin{eqnarray*}
&&A(\alpha,u,q):=x^{(\alpha)}x^{u}-\sum _{i\in
I(\alpha,u)}\Gamma_{i}^{q}(x^{(\alpha)}x^{u}),
\\
&& B(\alpha,u,\lambda,q):=(-1)^{|
u|}(n\lambda-\mathrm{zd}(x^{(\alpha)}x^{u}))\nabla_{q}(x^{(\alpha)}x^{u}),
\\
&&
E(\alpha,u):= D_{KO}(x^{(\alpha)}x^{u}),\\
&&
 E(\alpha,u,q):=D_{KO}\big(A(\alpha,u,q)\big),\\
&&
 E(\alpha,u,2n+1):= D_{KO}(x^{(\alpha)}x^{u}x_{2n+1}),\\
&&
 E(\alpha,u,2n+1,q):=D_{KO}\big(x^{(\alpha)}x^{u}x_{2n+1}+B(\alpha,u,\lambda,q)\big),\\
&&
 G(\alpha,u,2n+1,q):=D_{KO}\big(A(\alpha,u,q)x_{2n+1}+B(\alpha,u,\lambda,q)\big)
 \end{eqnarray*}
 and
\begin{eqnarray*}
&&S_{1}:=\{  E(\alpha,u)\mid  I(\alpha,u)=\emptyset,
(\alpha,u)\not=(0,0)\},
\\
&& S_{2}:=\{ E(\alpha,u,q)\mid x^{(\alpha)}x^{u}\in
\mathcal{D}^{\ast}, q\in \widetilde{I}(\alpha,u)\},
\\
&& S_{3}:=\{ E(\alpha,u,2n+1,q(\alpha,u))\mid  I(\alpha,u)=\emptyset
, \widetilde{I}(\alpha,u)\neq\emptyset\},
\\
&& S_{4}:=\{ G(\alpha,u,2n+1,q)\mid x^{(\alpha)}x^{u}\in
\mathcal{D}^{\ast}, q\in \widetilde{I}(\alpha,u)\},
\\
&& S_{5}:=\{ E(\alpha,u,2n+1)\mid I(\alpha,u)=\emptyset,
\widetilde{I}(\alpha,u)=\emptyset,
n\lambda-\mathrm{zd}(x^{(\alpha)}x^{u})=0\; \mbox{in}\;
\mathbb{F}\}.
\end{eqnarray*}
One can easily verify that $ \cup_{i=1}^{5} S_{i} \subset
\mathfrak{g}''$ and
$\mathrm{span}_{\mathbb{F}}S_{i}\cap\mathrm{span}_{\mathbb{F}}(\cup_{j\not=i}
S_{j})=0. $

In the following, we use frequently a decomposition: Given $j\in \overline{n+1,2n+ 1}$ any element  $f\in \mathcal{O}$ can be uniquely  written as
\begin{equation*}
f=f_{0}x_{j}+f_{1}\quad \mbox{where}\quad
 \partial_{j}(f_{0})=\partial_{j}(f_{1})=0,
 \end{equation*}
called the $x_{j}$-\textit{decomposition} of $f.$ The $x_{2n+1}$-decomposition of $f\in \mathcal{O}$ is
\begin{equation}
\label{decomposition} f=f_{0}x_{2n+1}+f_{1}\quad \mbox{where}\quad
 \partial_{2n+1}(f_{0})=\partial_{2n+1}(f_{1})=0.
\end{equation}
\begin{theorem}\label{t5}
$\mathfrak{g}''$ is spanned by
 $\{D_{KO}(1)\}$ and $\cup_{i=1}^{5} S_{i}.$
\end{theorem}
\begin{proof}
Let $D_{KO}(f)$ be an arbitrary element of $\mathfrak{g}'',$ where
$f\in\mathcal{O}$ such that $\mathrm{div}_{\lambda}(f)=0.$ We want
to show that $D_{KO}(f)$ is a linear combination of $\{D_{KO}(1)\}$
and $\cup_{i=1}^{5} S_{i}.$ Consider the $x_{2n+1}$-decomposition
(\ref{decomposition})  of $f$. We have
$$
2(-1)^{\mathrm{p}(f)}\big(\Delta(f_{0})x_{2n+1}+\Delta(f_{1})+(-1)^{\mathrm{p}(f_{0})}(\mathfrak{D}-n\lambda\mathrm{id}_{\mathcal{O}})(f_{0})\big)
=\mathrm{div}_{\lambda}(f)=0.
$$
It follows that
\begin{equation}\label{liue1021}
  \Delta(f_{0})=0
  \end{equation}
  and
\begin{equation}\label{liue1022}
  \Delta(f_{1})=(-1)^{\mathrm{p}(f_{0})}(n\lambda\mathrm{id}_{\mathcal{O}}-\mathfrak{D})(f_{0}).
  \end{equation}
From (\ref{liue1021}) one sees that $f_{0}\in \mathbb{F}$ or $D_{KO}(f_{0})\in
\mathcal{K}(n,n+1;\lambda,\underline{t}).$
 From (\ref{liue2}) we know that
 $\mathcal {K}(n,n+1;\lambda,\underline{t})\cong SHO'(n,n;\underline{t}).$ Then by \cite[Theorem 2.7]{lh},
$\mathcal {K}(n,n+1;\lambda,\underline{t})$ is spanned by $S_{1}\cup
S_{2}.$ Thus, we can suppose
\begin{equation}\label{liue3} f_{0}=\sum_{I(\alpha,u)=\emptyset}a_{(\alpha,u)}x^{(\alpha)}x^{u}+
\sum_{ x^{(\alpha)}x^{u}\in \mathcal{D}^{\ast}, q\in \widetilde{I}(\alpha,u)}
  a_{(\alpha,u)}A(\alpha,u,q),
\end{equation}
where $a_{(\alpha,u)}\in \mathbb{F}.$

Suppose  $n\lambda-\mathrm{zd}(x^{(\alpha)}x^{u})\neq 0$ in
$\mathbb{F}$ for some $\alpha, u$ with $I(\alpha,u)=\emptyset $ and
 $\widetilde{I}(\alpha,u)=\emptyset$.
It follows  from (\ref{liue1022})  that $a_{(\alpha,u)}=0.$   Thus (\ref{liue3}) can be rewritten as
 \begin{eqnarray}\label{liue4}
 f_{0}&=&\sum_{I(\alpha,u)=\emptyset, \widetilde{I}(\alpha,u)=\emptyset\atop
n\lambda-\mathrm{zd}(x^{(\alpha)}x^{u})=0}a_{(\alpha,u)}x^{(\alpha)}x^{u}+\sum_{I(\alpha,u)=\emptyset,
 \widetilde{I}(\alpha,u)\not=\emptyset}a_{(\alpha,u)}x^{(\alpha)}x^{u}\nonumber\\
 &&+\sum_{ x^{(\alpha)}x^{u}\in \mathcal{D}^{\ast}, q\in \widetilde{I}(\alpha,u)}
  a_{(\alpha,u)}A(\alpha,u,q).
  \end{eqnarray}
Let
\begin{eqnarray}\label{liue5}
g&:=&f_{1}-\sum_{I(\alpha,u)=\emptyset, \widetilde{I}(\alpha,u)\not=\emptyset}a_{(\alpha,u)}B(\alpha,u,\lambda,q(\alpha,u))\nonumber\\
&&-\sum_{ x^{(\alpha)}x^{u}\in \mathcal{D}^{\ast}, q\in
\widetilde{I}(\alpha,u)}
  a_{(\alpha,u)}B(\alpha,u,\lambda,q).
  \end{eqnarray}
  Then
\begin{eqnarray*}
\Delta(f_{1})
&=&\Delta(g)
 +\Delta\Big(\sum_{I(\alpha,u)=\emptyset, \widetilde{I}(\alpha,u)\not=\emptyset}a_{(\alpha,u)}B(\alpha,u,\lambda,q(\alpha,u))\\
&&+\sum_{ x^{(\alpha)}x^{u}\in \mathcal{D}^{\ast}, q\in
\widetilde{I}(\alpha,u)}
  a_{(\alpha,u)}B(\alpha,u,\lambda,q)\Big)\\
&=&\Delta(g)
+(-1)^{\mathrm{p}(f_{0})}(n\lambda\mathrm{id}_{\mathcal{O}}-\mathfrak{D})(f_{0}).
\end{eqnarray*}
Then by (\ref{liue1022}), we have $\Delta(g)=0.$ Hence $g\in \mathbb{F}$ or
$D_{KO}(g)\in \mathcal {K}(n,n+1;\lambda,\underline{t}).$
Consequently,
\begin{eqnarray}\label{liue6}
D_{KO}(g)\in \mathrm{span}_{\mathbb{F}}(S_{1}\cup
S_{2}\cup\{D_{KO}(1)\}).
 \end{eqnarray}
From (\ref{liue4}) and (\ref{liue5}) it follows that
\begin{eqnarray*}
f&=&f_{0}x_{2n+1}+f_{1}\\
&=&f_{0}x_{2n+1}+\sum_{I(\alpha,u)=\emptyset, \widetilde{I}(\alpha,u)\not=\emptyset}a_{(\alpha,u)}B(\alpha,u,\lambda,q(\alpha,u))\\
&&+\sum_{ x^{(\alpha)}x^{u}\in \mathcal{D}^{\ast}, q\in
\widetilde{I}(\alpha,u)}
  a_{(\alpha,u)}B(\alpha,u,\lambda,q)+g\\
&=&\sum_{I(\alpha,u)=\emptyset,
\widetilde{I}(\alpha,u)=\emptyset\atop
n\lambda-\mathrm{zd}(x^{(\alpha)}x^{u})=0}a_{(\alpha,u)}x^{(\alpha)}x^{u}x_{2n+1}+\sum_{I(\alpha,u)
=\emptyset, \widetilde{I}(\alpha,u)\not=\emptyset}a_{(\alpha,u)}(x^{(\alpha)}x^{u}x_{2n+1}\\
&&+B(\alpha,u,\lambda,q(\alpha,u)))+\sum_{ x^{(\alpha)}x^{u}\in
\mathcal{D}^{\ast}, q\in \widetilde{I}(\alpha,u)}
  a_{(\alpha,u)}(A(\alpha,u,q)x_{2n+1}\\
  &&+B(\alpha,u,\lambda,q))+g\\
  &=&\sum_{I(\alpha,u)=\emptyset,
\widetilde{I}(\alpha,u)=\emptyset\atop
n\lambda-\mathrm{zd}(x^{(\alpha)}x^{u})=0}a_{(\alpha,u)}E(\alpha,u,2n+1)+\sum_{I(\alpha,u)
=\emptyset, \widetilde{I}(\alpha,u)\not=\emptyset}a_{(\alpha,u)}E(\alpha,u,2n+1,q(\alpha,u))\\
&&+\sum_{ x^{(\alpha)}x^{u}\in \mathcal{D}^{\ast}, q\in
\widetilde{I}(\alpha,u)}
  a_{(\alpha,u)}G(\alpha,u,2n+1,q)+g.
\end{eqnarray*}
This combining  (\ref{liue6}) shows that

$$D_{KO}(f)\, \mbox{is a linear combination of}\, \cup_{i=1}^{5} S_{i}\cup\{D_{KO}(1)\}.$$
 \end{proof}
\begin{lemma}\label{l6}
Suppose  $ D_{KO}(f)\in \mathfrak{g}''$  and $f=f_{0}x_{2n+1}+f_{1}$
is the $x_{2n+1}$-decomposition. Then $[D_{KO}(f),D_{KO}(1)]=2
D_{KO}(f_{0}).$
\end{lemma}
\begin{proof} Using (\ref{liuee1}), one may directly compute.
\end{proof}

\begin{lemma}\label{l8}
If  $ D_{KO}(f)\in [\mathcal {K}(n,n+1;\lambda,\underline{t}),
\mathcal {K}(n,n+1;\lambda,\underline{t})],$  then $f$ does not contain any nonzero
monomials of the form $x^{(\alpha)}x^{u}$ such that $I(\alpha,u)=\emptyset,$
$\widetilde{I}(\alpha,u)=\emptyset.$
\end{lemma}
\begin{proof} Write $D_{KO}(f)=[D_{KO}(g),D_{KO}(h)] $, where $g, h$ satisfy the conditions  (\ref{liue1}).
Then by (\ref{liue2}),
$$
T_{H}(f)=[T_{H}(g), T_{H}(h)]
$$
 and the conclusion follows from \cite[Proposition 3.4]{lh}.
\end{proof}
Put
\begin{eqnarray*}
&&\mathfrak{A}_{1}:=\{x^{(\alpha)}x^{u}\mid
I(\alpha,u)=\widetilde{I}(\alpha,u)=\emptyset\},
\\
&&\mathfrak{A}_{2}:=\{x^{(\alpha)}x^{u}\mid I(\alpha,u)=\emptyset,
\widetilde{I}(\alpha,u)\neq\emptyset\}.
\end{eqnarray*}
Note that
\begin{eqnarray}\label{e13}
x^{(\alpha)}x^{u}\in \mathfrak{A}_{1}\Longleftrightarrow
\alpha_{i}=\pi_{i} \;\mbox{or}\; 0 \;\mbox{for all}\; i\in
\overline{1,n} \;\mbox{and}\; u=\{j'\mid \alpha_{j}=0, j\in
\overline{1,n}\}.
\end{eqnarray}
It is clear that
$$
S_{1}\cup\{D_{KO}(1)\}=\{D_{KO}(f)\mid f\in \mathfrak{A}_{1}\}\cup
\{D_{KO}(f)\mid f\in \mathfrak{A}_{2}\}.
$$
Given $r\in \overline{1,n},$ let
$$
\mathbf{J}(0):=\emptyset,
$$
$$
\mathbf{J}(r):=\{(i_{1}, \ldots, i_{r})\mid  1\leq i_{1}<\cdots< i_{r}  \leq n\}.
$$
For $(i_{1}, \ldots, i_{r})\in \mathbf{J}(r)$, $r\in
\overline{0,n},$ let
$$
X(i_{1}, \ldots, i_{r}):=x^{(\pi_{i_{1}}\varepsilon_{i_{1}}+\cdots
 +\pi_{i_{r}}\varepsilon_{i_{r}})}x^{\langle 1',\ldots, n'\rangle-\langle i'_{1}, \ldots, i'_{r}\rangle}.
$$
We list some  technical formulas, which will be used later. Suppose
$I(\alpha^{1},u^{1})=\emptyset,$ $I(\alpha^{2},u^{2})=\emptyset,$
$\alpha^{1}+\alpha^{2}=\alpha,$ $u^{1}+ u^{2}=u,$ $q\in
\widetilde{I}(\alpha,u).$ Then
\begin{eqnarray}\label{e2.15}
[E(\alpha^{1},u^{1},2n+1,q),
     E(\alpha^{2},u^{2},2n+1,q)]= \left\{\begin{array}{ll}
\gamma E(\alpha,u,2n+1,q) &\mbox{if}\; I(\alpha,u)=\emptyset\\
\gamma G(\alpha,u,2n+1,q) &\mbox{if}\; I(\alpha,u)\neq\emptyset,
\end{array}\right.
\end{eqnarray}
where
\begin{eqnarray*}\gamma:&=&\pm\Big[(n\lambda-\mathrm{zd}(\alpha^{2},u^{2})){\alpha\choose\alpha^1
-\varepsilon_{q}}+(\mathrm{zd}(\alpha^{1},u^{1})-\mathrm{zd}(\alpha^{2},u^{2}))
{\alpha\choose\alpha^1}\\
&&-(n\lambda-\mathrm{zd}(\alpha^{1},u^{1})){\alpha\choose\alpha^1+\varepsilon_{q}}\Big].
\end{eqnarray*}

Suppose $f\in \mathcal {O}(n,n;\underline{t})$ is $q$-integral,
where $q \in \overline{1,n}.$ Put
$$Y(f,q):=D_{KO}(fx_{2n+1}+(-1)^{\mathrm{p}(f_{0})}(n\lambda-\mathrm{zd}(f_{0}))\nabla_{q}(f_{0})).$$
In the sequel, once the symbol $Y(f,q) $ appears, it will
  impliy that $f\in \mathcal {O}(n,n;\underline{t})$ is $q$-integral.

NOTICE that in  Lemmas \ref{l10}, \ref{l1},  Theorems  \ref{t9} and
\ref{t8} we will always  assume that $D_{KO}(f_{0}),$ $
D_{KO}(g_{0})\in S_{1}\cup S_{2}\cup\{D_{KO}(1)\}.$
\begin{lemma}\label{l10}
We have
\begin{eqnarray}\label{e2.1}
&&[Y(f_{0},q),D_{KO}(g_{0})]\nonumber\\
&=&D_{KO}\Big(\sum_{i=1}^{2n}(-1)^{\mu(i')\mathrm{p}(f_{0}x_{2n+1})}\partial_{i'}
(f_{0}x_{2n+1})\partial_{i}(g_{0})+(2-\mathrm{zd}(g_{0}))
f_{0}g_{0}\nonumber\\
&& +(-1)^{\mathrm{p}(f_{0})}(n\lambda-\mathrm{zd}(f_{0}))
\sum_{i=1}^{2n}(-1)^{\mu(i')\mathrm{p}(\nabla_{q}(f_{0}))}\partial_{i'}
(\nabla_{q}(f_{0}))\partial_{i}(g_{0})\Big)
\end{eqnarray}
and
\begin{eqnarray}\label{e2.2}
&&[Y(f_{0},q),
Y(g_{0},r)]\nonumber\\
&=&D_{KO}\Big((-1)^{\mathrm{p}(f_{0})+\mathrm{p}(g_{0})}(n\lambda-
\mathrm{zd}(f_{0}))(\mathrm{zd}(\nabla_{q}(f_{0}))-2)\nabla_{q}(f_{0})g_{0}\nonumber\\&&+(-1)^{\mathrm{p}(f_{0})
+\mathrm{p}(g_{0})}
(n\lambda-\mathrm{zd}(f_{0}))(n\lambda-\mathrm{zd}(g_{0}))
\sum_{i=1}^{2n}(-1)^{\mu(i)\mathrm{p}(\nabla_{q}(f_{0}))}
\partial_{i'}(\nabla_{q}(f_{0}))\partial_{i}(\nabla_{r}(g_{0}))\nonumber
\\
&&+(-1)^{\mathrm{p}(g_{0})}(n\lambda-\mathrm{zd}(g_{0}))\sum_{i=1}^{2n}(-1)^{\mu(i')
\mathrm{p}(f_{0}x_{2n+1})}
\partial_{i'}(f_{0}x_{2n+1})\partial_{i}(\nabla_{r}(g_{0}))\nonumber
\\
&&+(-1)^{\mathrm{p}(f_{0})}
(n\lambda-\mathrm{zd}(f_{0}))\sum_{i=1}^{2n}
(-1)^{\mu(i')\mathrm{p}(\nabla_{q}(f_{0}))}\partial_{i'}(\nabla_{q}(f_{0}))\partial_{i}(g_{0}x_{2n+1})\nonumber
\\
&&+(-1)^{\mathrm{p}(g_{0})}(n\lambda-\mathrm{zd}(g_{0}))(2+\mathrm{zd}(\nabla_{r}(g_{0})))f_{0}\nabla_{r}(g_{0})\nonumber
\\
&&+(\mathrm{zd}(f_{0})-\mathrm{zd}(g_{0}))f_{0}g_{0}x_{2n+1}\Big).
\end{eqnarray}
\end{lemma}
\begin{proof}
Using (\ref{liuee1}), one can directly compute.
\end{proof}
Given  $\lambda\in \mathbb{F}$  and $l\in \mathbb{Z},$ put
$$
\mathfrak{S}_{l}(\lambda, n):=\{k\in \overline{0,n}\mid
n\lambda-n+2k+l=0 \in\mathbb{F}\},
$$
the set of all the integer solutions between $0$ and $n$ of the
equation that $n\lambda-n+2x+l=0$ in $\mathbb{F}.$ Put
$$G:=G(\pi-\varepsilon_{1},\langle 2',\ldots,n'\rangle,2n+1,1).$$
\begin{lemma}\label{l1}
If $ [Y(f_{0},l),Y(g_{0},r)]=sG $ for some $0\neq s\in \mathbb{F},$
then there is $k\in \overline{1,n}$ such that $f_{0}$ and $g_{0}$
are $k$-integral.
\end{lemma}
\begin{proof}
If there exists no such $k$, then
$$
\nabla_{k}(f_{0}g_{0})=\nabla_{k}(\partial_{i'}(f_{0})\partial_{i}(\nabla_{r}(g_{0})))=\nabla_{k}(\partial_{i'}(\nabla_{l}(f_{0}))\partial_{i}(g_{0}))=0
$$
for all $k\in \overline{1,n}.$ Since $ [Y(f_{0},l),Y(g_{0},r)]=sG, $
it follows that
$$f_{0}g_{0}=\partial_{i'}(f_{0})\partial_{i}(\nabla_{r}(g_{0}))=\partial_{i'}(\nabla_{l}(f_{0}))\partial_{i}(g_{0})=0.$$
Then by (\ref{e2.2})  we have  $s=0.$ This contradicts the assumption
that $s\neq0.$ The proof is complete.
\end{proof}
\begin{theorem}\label{t9} If $n\lambda+1\not=0$ in $\mathbb{F}$ or
$\mathfrak{S}_{0}(n,\lambda)\neq\emptyset,$ then
$$
\mathfrak{g}''=\mathfrak{g}'\oplus\mathrm{span}_{\mathbb{F}}S_{5}\\
\oplus\sum_{r\in\mathfrak{S}_{2}(\lambda, n)\atop
(i_{1},\ldots,i_{r})\in \mathbf{J}(r)}\mathbb{F} D_{KO}(X(i_{1},
\ldots, i_{r})).
$$
If $n\lambda+1=0$ in $\mathbb{F}$ and
$\mathfrak{S}_{0}(n,\lambda)=\emptyset,$ then
$$
\mathfrak{g}''=\mathfrak{g}'\oplus\mathrm{span}_{\mathbb{F}}S_{5}\\
\oplus\sum_{r\in\mathfrak{S}_{2}(\lambda, n)\atop
(i_{1},\ldots,i_{r})\in \mathbf{J}(r)}\mathbb{F} D_{KO}(X(i_{1},
\ldots, i_{r}))\oplus\mathbb{F}G.
$$
\end{theorem}
\begin{proof}
In the light of  Theorem \ref{t5}, our discussion is divided into
six parts.

\textit{Part 1}.  Assert that $\{D_{KO}(f)\mid f\in \mathfrak{A}_{2}\}\cup
S_{2}\subseteq\mathfrak{g}'.$ This follows from Lemma \ref{l6}.

\textit{Part 2}. Assert that $S_{3} \subseteq \mathfrak{g}'.$ Consider the
elements in $S_{3}$, say, $E(\alpha,u,2n+1,q(\alpha,u)).$ Take
$q\neq q_{(\alpha,u)}'.$ Since
\begin{eqnarray*}
E(\alpha,u,2n+1,q_{(\alpha,u)}) \equiv
[E(\alpha+\varepsilon_{q_{(\alpha,u)}},u),
E(0,q_{(\alpha,u)}',2n+1,q)] \pmod{S_{2}},
\end{eqnarray*}
we have $E(\alpha,u,2n+1,q_{(\alpha,u)})\in \mathfrak{g}'.$ Thus,
$S_{3}\subseteq \mathfrak{g}'.$

\textit{Part 3}. Assert that $S_{4}\setminus \{\pm
G\}\subseteq\mathrm{alg}_{\mathbb{F}}(S_{3}\cup
\{D_{KO}(1)\})\subseteq \mathfrak{g}'.$ Consider
$G(\alpha,u,2n+1,q)\in S_{4}$. Suppose $I(\alpha,u)=\{i_{1},\ldots
,i_{k}\}\neq \emptyset$ and $u=\{i_{1}',\ldots ,i_{k}',
i_{k+1}',\ldots, i_{r}'\}\neq \emptyset.$ Write
\begin{eqnarray*}
&&(\alpha^{1},u^{1}):=
(\alpha_{i_{1}}\varepsilon_{i_{1}}+\cdots+\alpha_{i_{k}}\varepsilon_{i_{k}},\langle
i_{k+1}',\ldots ,i_{r}'\rangle);\\
&&
 (\alpha^{2},u^{2}):=(\alpha-\alpha^{1},u-u^{1});\\
&&(\alpha^{3},u^{3}):=
(\alpha_{i_{2}}\varepsilon_{i_{2}}+\cdots+\alpha_{i_{k}}\varepsilon_{i_{k}},\langle
i_{1}',i_{k+1}',\ldots ,i_{r}'\rangle);\\
&&
 (\alpha^{4},u^{4}):=(\alpha-\alpha^{3},u-u^{3}).
\end{eqnarray*}
By (\ref{e2.15}),  we have
\begin{eqnarray}\label{e2.17}
&&[E(\alpha^{1},u^{1},2n+1,q), E(\alpha^{2},u^{2},2n+1,q)]\nonumber\\
&=& \big(n\lambda
\alpha_{q}-(\alpha_{q}+1)\mathrm{zd}(x^{(\alpha^{2})}x^{u^{2}})+\mathrm{zd}(x^{(\alpha^{1})}x^{u^{1}})\big)
G(\alpha,u,2n+1,q)
\end{eqnarray}
and
\begin{eqnarray}\label{e2.18}
&&[E(\alpha^{3},u^{3},2n+1,q), E(\alpha^{4},u^{4},2n+1,q)]\nonumber\\
&=&\big(n\lambda
\alpha_{q}-(\alpha_{q}+1)\mathrm{zd}(x^{(\alpha^{4})}x^{u^{4}})+\mathrm{zd}(x^{(\alpha^{3})}x^{u^{3}})\big)
G(\alpha,u,2n+1,q).
\end{eqnarray}
Note that
\begin{eqnarray}\label{e2.19}
 (\alpha_{q}+2)(\alpha_{i_{1}}-1)&=&\big(n\lambda
\alpha_{q}-(\alpha_{q}+1)\mathrm{zd}(x^{(\alpha^{2})}x^{u^{2}})+\mathrm{zd}(x^{(\alpha^{1})}x^{u^{1}})\big)\nonumber
\\
&&-\big(n\lambda
\alpha_{q}-(\alpha_{q}+1)\mathrm{zd}(x^{(\alpha^{4})}x^{u^{4}})+\mathrm{zd}(x^{(\alpha^{3})}x^{u^{3}})\big).
\end{eqnarray}

Let us show the assertion in Part 3.

\textit{Subpart 3.1}. Assert that $ G(\alpha,u,2n+1,q)\in
\mathrm{alg}_{\mathbb{F}}(S_{3}\cup \{D_{KO}(1)\}) $ if
$\alpha_{q}\not \equiv-2\pmod{p}$ and $\alpha_{i_{1}}\not \equiv
1\pmod{p}.$ From (\ref{e2.19}) we know that
$$n\lambda
\alpha_{q}-(\alpha_{q}+1)\mathrm{zd}(x^{(\alpha^{2})}x^{u^{2}})+\mathrm{zd}(x^{(\alpha^{1})}x^{u^{1}})
$$
and
$$n\lambda
\alpha_{q}-(\alpha_{q}+1)\mathrm{zd}(x^{(\alpha^{4})}x^{u^{4}})+\mathrm{zd}(x^{(\alpha^{3})}x^{u^{3}})
$$
cannot be all zero in $\mathbb{F}$. By (\ref{e2.17}) and (\ref{e2.18}),
we have
$$ G(\alpha,u,2n+1,q)\in \mathrm{alg}_{\mathbb{F}}(S_{3}\cup
\{D_{KO}(1)\}).
$$

\textit{Subpart 3.2}. If $\alpha_{q}\not \equiv-2\pmod{p}$ and
$\alpha_{i_{1}}\equiv 1\pmod{p},$ then by Subpart 3.1,
$$
G(\alpha+\varepsilon_{i_{1}},u,2n+1,q)\in
\mathrm{alg}_{\mathbb{F}}(S_{3}\cup \{D_{KO}(1)\}).
 $$ Since
$E(0,i_{1}')\in \mathrm{alg}_{\mathbb{F}}(S_{3}\cup \{D_{KO}(1)\}),$
we have
$$
G(\alpha,u,2n+1,q)=[G(\alpha+\varepsilon_{i_{1}},u,2n+1,q),E(0,i_{1}')]\in
\mathrm{alg}_{\mathbb{F}}(S_{3}\cup \{D_{KO}(1)\}).
$$

\textit{Subpart 3.3}. If $\alpha_{q}\equiv-2\pmod{p}$ and
$-2n\lambda+\mathrm{zd}(\alpha,u)\not=0$ in $\mathbb{F},$ then by
(\ref{e2.15}) , we have $$ G(\alpha,u,2n+1,q)\in
\mathrm{alg}_{\mathbb{F}}(S_{3}\cup \{D_{KO}(1)\}).$$

\textit{Subpart 3.4}. If
 $\alpha_{q}\equiv-2\pmod{p}$ and
$-2n\lambda+\mathrm{zd}(\alpha,u)=0$ in $\mathbb{F}$ and $
G(\alpha,u,2n+1,q)\neq\pm G, $ then there is $i\in \overline{1,n}$
such that
$$
G(\alpha+\varepsilon_{i},u,2n+1,q)\in
\mathrm{alg}_{\mathbb{F}}(S_{3}\cup \{D_{KO}(1)\})$$
or
$$G(\alpha,u+\langle i'\rangle,2n+1,q)\in
\mathrm{alg}_{\mathbb{F}}(S_{3}\cup \{D_{KO}(1)\}).
$$
Since $E(\varepsilon_{i},0), E(0,i')\in
\mathrm{alg}_{\mathbb{F}}(S_{3}\cup \{D_{KO}(1)\}),$ we have
$$
G(\alpha,u,2n+1,q)\equiv[G(\alpha+\varepsilon_{i},u,2n+1,q),E(0,i')]\pmod
{\mathrm{alg}_{\mathbb{F}}(S_{3}\cup \{D_{KO}(1)\})}
$$
or
$$
G(\alpha,u,2n+1,q)\equiv[G(\alpha,u+\langle
i'\rangle,2n+1,q),E(i,0)]\pmod {\mathrm{alg}_{\mathbb{F}}(S_{3}\cup
\{D_{KO}(1)\})}.
$$
Therefore, $G(\alpha,u,2n+1,q)\in
\mathrm{alg}_{\mathbb{F}}(S_{3}\cup \{D_{KO}(1)\}).$
It follows from  Subparts 3.1--3.4 that $S_{4}\setminus \{\pm
G\}\subseteq \mathrm{alg}_{\mathbb{F}}(S_{3}\cup \{D_{KO}(1)\}).$

\textit{Part 4}. Assert that $  G\in \mathfrak{g}'\Longleftrightarrow
n\lambda+1\not=0\;\mbox{in}\;\mathbb{F}$ or
$\mathfrak{S}_{0}(n,\lambda)\neq\emptyset.$ Note that
\begin{eqnarray}\label{l11}
-2(n\lambda +1)G=[E(\pi-\varepsilon_{1},0,2n+1,1),E(0,\langle
2',\ldots,n'\rangle,2n+1,1)]\in \mathfrak{g}'.
\end{eqnarray}
If $n\lambda+1\not=0 $ in $\mathbb{F},$ then $ G\in \mathfrak{g}'.$ If
$\mathfrak{S}_{0}(n,\lambda)\neq\emptyset,$ then for $r\in
\mathfrak{S}_{0}(n,\lambda) $ we have
$$E(\pi_{1}+\cdots+\pi_{r},\langle (r+1)',\ldots,n',\rangle,2n+1)\in
\mathfrak{g''}.$$ Then
\begin{eqnarray*}
&&(-1)^{(n-r+1)(r-1)}G\\
&=&[E(\pi_{1}+\cdots+\pi_{r},\langle
(r+1)',\ldots,n',\rangle,2n+1),E(\pi_{r+1}+\cdots+\pi_{n},\langle
1',\ldots,r',\rangle)],
\end{eqnarray*}
that is, $G\in \mathfrak{g}'.$ Conversely, we consider two cases
separately:

 \textit{Case 1}. If $$
[Y(f_{0},q),D_{KO}(g_{0})]=sG, $$ then
$$[D_{KO}(f_{0}),D_{KO}(g_{0})]=sE(\pi-\varepsilon_{1},\langle 2',\ldots,n'\rangle,1),$$ where $s\in\mathbb{F}.$

\textit{Subcase 1.1}. Suppose $f_{0}=D_{KO}(x^{(\alpha)}x^{u})$,
$g_{0}=D_{KO}(x^{(\beta)}x^{v})$ with
 $I(\alpha,u)=I(\beta,v)=\emptyset$ and $\widetilde{I}(\alpha,u)\not=\emptyset,$
 $\widetilde{I}(\beta,v)=\emptyset.$ Note that
 $$[D_{KO}(f_{0}),D_{KO}(g_{0})]=\sum_{k=1}^{2n}(-1)^{\mu(k)|u|}\partial_{k}(x^{(\alpha)}x^{u})\partial_{k'}(x^{(\beta)}x^{v}). $$
 Assume that $\partial_{i}(x^{(\alpha)}x^{u})\partial_{i'}(x^{(\beta)}x^{v})\not=0$ for some $i\in \overline{1,2n}.$ By symmetry one may
 assume that $i\in \overline{1,n}.$ Then it is easily seen
 that $\partial_{i}(x^{(\alpha)}x^{u})\partial_{i'}(x^{(\beta)}x^{v})$ is of the form
 $sx^{(\pi-\varepsilon_{i})}x^{\omega-\langle i'\rangle}.$ It follows that
 $\alpha+\beta=\pi$ and $u+v=\omega$ (meaning  $u\cap v=\emptyset$). Noticing that $I(\alpha,u)=I(\beta,v)=\emptyset,$
one may easily deduce that $\widetilde{I}(\alpha,u)=\emptyset,$
contradicting our assumption. This shows that $s=0.$

\textit{Subcase 1.2}. Suppose
$f_{0}=D_{KO}(x^{(\alpha)}x^{u}-\sum_{i\in
I(\alpha,u)}\Gamma_{i}^{q}(x^{(\alpha)}x^{u})),$
$g_{0}=D_{KO}(x^{(\beta)}x^{v})$ with
 $I(\alpha,u)\neq\emptyset,$  $q\in \widetilde{I}(\alpha,u)\neq\emptyset,$
$I(\beta,v)=\emptyset$ and $\widetilde{I}(\beta,v)=\emptyset.$
Assume that $s\neq0.$ In the light of \cite [Lemma 3.2]{lh}, there
is $k\in \overline{1,2n}$ such that
$\partial_{k}(x^{(\alpha)}x^{u})\partial_{k'}(x^{(\beta)}x^{v})\neq0.$
Thus one may assume that
$\partial_{k}(x^{(\alpha)}x^{u})\partial_{k'}(x^{(\beta)}x^{v})$ is
of the form $sx^{(\pi-\varepsilon_{k})}x^{\omega-\langle
k'\rangle}.$ Consequently, we have
\begin{equation}\label{hee2.23}
\alpha+\beta=\pi \quad\mbox{and}\quad u+v=\omega.
\end{equation}
Note that $I(\beta,v)=\emptyset$ and
$\widetilde{I}(\beta,v)=\emptyset.$ By (\ref{e13}), we know there is
$(i_{1}, \ldots, i_{r})\in \mathbf{J}(r)$ such that
$x^{(\beta)}x^{v}=X(i_{1},\ldots,i_{r}).$ It follows from
(\ref{hee2.23}) that $x^{(\alpha)}x^{u}=X(i_{r+1},\ldots,i_{n}).$
This contradicts the assumption that $I(\alpha,u)\neq\emptyset,$
$\widetilde{I}(\alpha,u)\neq\emptyset.$

Now, in combination with \cite [Theorem 3.8]{lh}, one may prove $s=0.$

\textit{Case 2}.  Suppose
$f_{0}=x^{(\alpha)}x^{u}-\sum_{i\in
I(\alpha,u)}\Gamma_{i}^{l}(x^{(\alpha)}x^{u}),$
$g_{0}=x^{(\beta)}x^{v}-\sum_{i\in
I(\beta,v)}\Gamma_{i}^{r}(x^{(\beta)}x^{v}),$ where $l\in
\widetilde{I}(\alpha,u)\neq \emptyset,$
$r\in\widetilde{I}(\beta,v)\neq \emptyset.$  When
$I(\alpha,u)=I(\beta,v)=\emptyset,$  one sees that $f_{0}=x^{(\alpha)}x^{u} $ and
$g_{0}=x^{(\beta)}x^{v}.$  If
 $$
[Y(f_{0},l),Y(g_{0},r)]=sG,
$$
where $s\in \mathbb{F},$ by Lemma \ref{l1}, there is $k\in
\overline{1,n}$ such that $f_{0}$ and $g_{0}$ be $k$-integral. By
Case 1, we have
\begin{eqnarray*}
&&[Y(f_{0},l),Y(g_{0},r)]\\
&=&[Y(f_{0},k)+(-1)^{\mathrm{p}(f_{0})}(n\lambda-\mathrm{zd}(f_{0}))D_{KO}(\nabla_{l}(f_{0})-\nabla_{k}(f_{0})),\\
&&Y(g_{0},k)+(-1)^{\mathrm{p}(g_{0})}(n\lambda-\mathrm{zd}(g_{0}))D_{KO}(\nabla_{r}(g_{0})-\nabla_{k}(g_{0}))]\\
&=&[Y(f_{0},k),Y(g_{0},k)]\\
&=&sG=\pm sG(\pi-\varepsilon_{k},\langle 1',\ldots,n'\rangle-\langle
k'\rangle,2n+1,k).
\end{eqnarray*}
 Pay attention to the monomials in (\ref{e2.2})
which do not contain $x_{2n+1}$. Note that
$\mathrm{zd}(\alpha,u)+\mathrm{zd}(\beta,v)\equiv-2\pmod{p}.$
We have
\begin{eqnarray*}
&\pm& s(n\lambda+2)x^{(\pi)}x^{\langle1',\ldots,n'\rangle}\\
&=&(-1)^{|u|+|v|}(n\lambda-\mathrm{zd}(\alpha,u))\mathrm{zd}(\beta,v)x_{k'}x^{(\alpha+\varepsilon_{k})}x^{u}x^{(\beta)}x^{v}\\
&&+(-1)^{|u|+|v|}(n\lambda-\mathrm{zd}(\alpha,u))(n\lambda-\mathrm{zd}(\beta,v))x_{k'}x^{(\alpha)}x^{u}x^{(\beta+\varepsilon_{k})}x^{v}\\
&&-(-1)^{|u|+|v|}(n\lambda-\mathrm{zd}(\alpha,u))(n\lambda-\mathrm{zd}(\beta,v))x_{k'}x^{(\alpha+\varepsilon_{k})}x^{u}x^{(\beta)}x^{v}\\
&&+(-1)^{|u|+|v|}(n\lambda-\mathrm{zd}(\beta,v))(4+\mathrm{zd}(\beta,v))x_{k'}x^{(\alpha)}x^{u}x^{(\beta+\varepsilon_{k})}x^{v}\\
&=&\pm
2(n\lambda+1)(n\lambda+2)x^{(\pi)}x^{\langle1',\ldots,n'\rangle}\;\mbox{in}\;\mathbb{F}.
\end{eqnarray*}
It follows that $s=\pm2(n\lambda+1)\;\mbox{in}\;\mathbb{F}.$

\textit{Part 5}. Assert that for $r\in \overline{0,n},$
\begin{equation}\label{e1061}
D_{KO}(X(i_{1}, \ldots, i_{r}))\in \mathfrak{g}' \Longleftrightarrow
r\notin \mathfrak{S}_{2}(\lambda,n).
\end{equation}

\textit{Subpart 5.1}. By Lemma \ref{l8}, if $D_{KO}(h)\in [\mathcal
{K}(n,n+1;\lambda,\underline{t}), \mathcal
{K}(n,n+1;\lambda,\underline{t})],$ where $h\in \mathcal{O}$, then
$h$ has no nonzero monomials $X(i_{1}, \ldots, i_{r}).$

\textit{Subpart 5.2}. Pay attention to the monomials in (\ref{e2.2})
which do not contain $x_{2n+1}$. One sees that if  $[
Y(f_{0},q),Y(g_{0},r)]$$= D_{KO}(h),$ where $h\in \mathcal{O},$ then
$h$ has no nonzero monomials $X(i_{1}, \ldots, i_{r}).$

\textit{Subpart 5.3}. If there exist $Y(f_{0},q), D_{KO}(g_{0})$
such that $[ Y(f_{0},q),D_{KO}(g_{0})]= D_{KO}(h),$ where $h\in
\mathcal{O}$ and $h$ has  nonzero monomials $X(i_{1},\ldots,i_{r}),$
by (\ref{e2.1})  we have $f_{0}g_{0}=aX(i_{1}, \ldots, i_{r}),$
where $0\neq a\in \mathbb{F},$ and then
    \begin{eqnarray*}
   [Y(f_{0},q),D_{KO}(g_{0})]
=a(n\lambda-n+2r+2)X(i_{1}, \ldots, i_{r}).
\end{eqnarray*}

The assertion (\ref{e1061}) follows from Subparts 5.1--5.3.

\textit{Part 6}. Assert that $
\mathfrak{g}'\cap\mathrm{span}_{\mathbb{F}}S_{5}=0.$

\textit{Subpart 6.1}. Suppose $[Y(f_{0},q),Y(g_{0},r)]= D_{KO}(h) $
and $h\in \mathcal{O}$
  has  a nonzero monomial of the form
$x^{(\pi_{i_{1}}\varepsilon_{i_{1}})}\cdots
x^{(\pi_{i_{r}}\varepsilon_{i_{r}})}x_{i_{r+1}'}x_{i_{r+2}'}\cdots
x_{i_{n}'}x_{2n+1}.$  Then by (\ref{e2.2}), $f_{0}, g_{0}$ must be
of the form
$$f_{0}=x^{(\alpha)}x^{u},\;
 g_{0}=x^{(\beta)}x^{\upsilon},$$
 where
 $\alpha+\beta=\pi_{i_{1}}\varepsilon_{i_{1}}+\cdots+\pi_{i_{r}}\varepsilon_{i_{r}},$
 $u+\upsilon=\langle i_{r+1}'\ldots,i_{n}'\rangle.$
In this case,
\begin{eqnarray*}
&&[Y(f_{0},q),Y(g_{0},r)]=0.
\end{eqnarray*}

\textit{Subpart 6.2}. By (\ref{e2.1}), if $[
Y(f_{0},q),D_{KO}(g_{0})]= D_{KO}(h)$ and $h\in \mathcal{O}$ has a
nonzero monomial of the form
$x^{(\pi_{i_{1}}\varepsilon_{i_{1}})}\cdots
x^{(\pi_{i_{r}}\varepsilon_{i_{r}})}x_{i_{r+1}'}x_{i_{r+2}'}\cdots
x_{i_{n}'}x_{2n+1},$ then there are $0\neq a\in \mathbb{F}$ and $
h_{1}\in \mathcal {O}(n,n;\underline{t})$ such that
$$
[D_{KO}(f_{0}),D_{KO}(g_{0})]=aD_{KO}(h_{1}+X(i_{1}, \ldots,
i_{r})).
$$
This contradicts Lemma \ref{l8} and therefore, $
\mathfrak{g}'\cap\mathrm{span}_{\mathbb{F}}S_{5}=0.$ The proof is
complete.
\end{proof}
Put
\begin{eqnarray*}
\delta'_{n\lambda,-1}:= \left\{\begin{array}{ll}
0 &\mbox{if}\; n\lambda+1\not=0\quad \mbox{in}\; \mathbb{F}\\
1 &\mbox{if}\; n\lambda+1=0\quad \mbox{in}\; \mathbb{F}.
\end{array}\right.\nonumber
\end{eqnarray*}
\begin{theorem}\label{t8}
$\mathfrak{g}'=\mathfrak{g}\oplus\mathbb{F}\delta'_{n\lambda,-1}G.$
\end{theorem}
\begin{proof}
(1) By Parts 1--3 in the proof of Theorem \ref{t9} we have
$$\{D_{KO}(f)\mid f\in \mathfrak{A}_{2}\}\cup S_{2}\cup
S_{3}\subseteq \mathfrak{g},\;S_{4}\setminus \{\pm G\}\subseteq
\mathfrak{g}.$$

(2) Since
\begin{eqnarray*}
&&(n\lambda-n+2r+2)D_{KO}(X(i_{1}, \ldots, i_{r}))\\
&=&[E(\pi_{i_{1}}\varepsilon_{i_{1}},0,2n+1,i_{2}),E(\pi_{i_{2}}\varepsilon_{i_{2}}+\cdots
+\pi_{i_{r}}\varepsilon_{i_{r}},\langle i'_{r+1},\ldots,i'_{n}\rangle)],
\end{eqnarray*}
we have $ D_{KO}(X(i_{1}, \ldots, i_{r}))\in \mathfrak{g} $ for all
$r\notin \mathfrak{S}_{2}(\lambda,n).$

(3) We propose to show that
$$
G\in \mathfrak{g} \Longleftrightarrow n\lambda+1\not=0\quad \mbox{in}\;\mathbb{F}.
$$
By (\ref{l11}), we have $G\in \mathfrak{g}$ when
$n\lambda+1\not=0 \;\mbox{in}\;\mathbb{F}.$ The converse follows from Part 4
in the proof of Theorem \ref{t9}.
\end{proof}
For distinct $i, j, k\in \overline{1,n}$ and $q\in \overline{1,n},$ we list some  technical
formulas, which will be used later:
\begin{eqnarray}
 && [E(k_{i}\varepsilon_{i},0,2n+1,q), E(k_{j}\varepsilon_{j},0,2n+1,q)]
    =(k_{i}-k_{j})E(k_{i}\varepsilon_{i}+k_{j}\varepsilon_{j},0,2n+1,q),~~~~~\label{e2.3}\\
&&[E(2\varepsilon_{k},0,2n+1,q),
 E(0,\langle i'\rangle,2n+1,q)]=E(2\varepsilon_{k},\langle
 i'\rangle,2n+1,q),\label{e2.7}
    \\
    &&[E(2\varepsilon_{k},\langle i'\rangle,2n+1,q),
   E(0,\langle k'\rangle,2n+1,q)]=-2G(2\varepsilon_{k},\langle
   k',i'\rangle,2n+1,q),\label{e2.8}\\
 &&[G(2\varepsilon_{k},\langle k',i'\rangle,2n+1,q), E(0,\langle j'\rangle)]=E(2\varepsilon_{k},\langle
 k',i',j'\rangle,q),\label{e2.9}
   \\
   &&[E(2\varepsilon_{k},\langle
  k',i',j'\rangle,q),E(j,0)]=-E(2\varepsilon_{k},\langle
  k',i'\rangle,q),\label{e2.10'}
  \\
  &&[E(2\varepsilon_{k},\langle
  k',i',j'\rangle,q),E(0,\langle k'\rangle)]=E(\varepsilon_{k},\langle
  k',i',j'\rangle,q),\label{e2.10}
   \\
   &&[E(2\varepsilon_{k},\langle i'\rangle,2n+1,q),
    E(0,\langle k'\rangle)]=
     E(\varepsilon_{k},\langle i'\rangle,2n+1,q)-E(2\varepsilon_{k},\langle
     k',i'\rangle,q),\label{e2.11}
    \\
    &&[E(\varepsilon_{k},\langle i'\rangle,2n+1,q),
     E(0,\langle j'\rangle)]=-E(\varepsilon_{k},\langle i',j'\rangle),\;\;q\neq
     j,\label{e2.12}
     \\
     &&
  [E(\varepsilon_{k},\langle i',j'\rangle),
   E(0,\langle k'\rangle,2n+1,q)]=
    E(0,\langle i',j'\rangle,2n+1,q)-E(\varepsilon_{k},\langle k',i',j'\rangle,q),\label{e2.14}
   \end{eqnarray}
   and
   \begin{eqnarray}
  &&[E(k_{i}\varepsilon_{i}+(k_{j}+1)\varepsilon_{j},0,2n+1,q), E(0,\langle j'\rangle,2n+1,q)]\nonumber
     \\
     &=&(k_{i}+k_{j})G(k_{i}\varepsilon_{i}+(k_{j}+1)\varepsilon_{j},\langle
     j'\rangle,2n+1,q),\label{e2.4}
     \\
   &&[E(\pi_{i}\varepsilon_{i}+(\pi_{j}-1)\varepsilon_{j},0,2n+1,q),
     E(\varepsilon_{j},0,2n+1,q)]
     \nonumber
     \\&=&\pi_{j}(\pi_{i}+\pi_{j}-2)E(\pi_{i}\varepsilon_{i}+\pi_{j}\varepsilon_{j},0,2n+1,q),\label{e2.5}
   \\
    &&[E(k_{i}\varepsilon_{i}+(k_{j}+1)\varepsilon_{j},0),
     E(0,\langle j'\rangle,2n+1,i)]
    \nonumber
     \\&=&E(k_{i}\varepsilon_{i}+k_{j}\varepsilon_{j},0,2n+1,j)
     +(k_{i}+1)(n\lambda-1)E((k_{i}+1)\varepsilon_{i}+k_{j}\varepsilon_{j},\langle
     i'\rangle,j).\label{e2.6}
     \end{eqnarray}
Let
\begin{eqnarray*}
&&T=\{E(k_{i}\varepsilon_{i},0,2n+1,q)\mid 1\leq i\leq n,  0\leq
k_{i}\leq p^{t_{i}}-1, q\in \widetilde{I}(k_{i}\varepsilon_{i},0)\};\\
&&S=\{E(0,i',2n+1,q)\mid
 1\leq i\leq n, q\in \widetilde{I}(0,\langle i'\rangle)\}.
\end{eqnarray*}
Note that $T,S\subseteq S_{3}.$
\begin{theorem}\label{t10}
$\mathfrak{g}$ is generated by $T\cup S\cup \{D_{KO}(1)\}.$
\end{theorem}
\begin{proof}
Let $Y:=\mathrm{alg}_{\mathbb{F}}(T\cup S \cup \{D_{KO}(1)\}).$ By
Theorem \ref{t9}, we may complete the proof by the following four steps.

\textit{Step 1}. Use induction on the number of
variables of
$E(\alpha,u,2n+1,q(\alpha,u))\in S_{3} $ to show that
$S_{3}\subseteq Y.$

\textit{Case 1.1.} Assert that for distinct $i, j,q\in
\overline{1,n},$
$$E(k_{i}\varepsilon_{i}+k_{j}\varepsilon_{j},0,2n+1,q)\in Y.$$
If $k_{i}-k_{j}\not \equiv 0\pmod{p},$ the conclusion follows from (\ref{e2.3}).
If $k_{i}-k_{j}\equiv 0\pmod{p},$  we consider two cases separately:

 \textit{Subcase 1.1.1}. Suppose $k_{i}\neq \pi_{i}$ or $k_{j}\neq \pi_{j},$ say,
$k_{j}\neq \pi_{j}.$ From the above, we have
$$E(k_{i}\varepsilon_{i}+(k_{j}+1)\varepsilon_{j},0,2n+1,q)\in Y.$$
By Lemma \ref{l6},
     $$
E(k_{i}\varepsilon_{i}+(k_{j}+1)\varepsilon_{j},0)\in Y.
     $$
If $k_{i}+k_{j}\not\equiv0 \pmod{p},$ then by (\ref{e2.4}) we have,
\begin{eqnarray*}
&&G(k_{i}\varepsilon_{i}+(k_{j}+1)\varepsilon_{j},\langle
j'\rangle,2n+1,q)\in Y.
\end{eqnarray*}
It follows from  Lemma \ref{l6} that
$$E(k_{i}\varepsilon_{i}+(k_{j}+1)\varepsilon_{j},\langle j'\rangle,q)\in
Y.$$
If $k_{i}+k_{j}\equiv 0 \pmod{p},$ then $k_{i}\leq \pi_{i}-2$ or
$k_{j}\leq \pi_{j}-2.$ Without loss of generality, suppose
$k_{i}\leq \pi_{i}-2.$ From the above, we have
$$E((k_{i}+1)\varepsilon_{i}+(k_{j}+1)\varepsilon_{j},\langle j'\rangle,q)\in
Y.$$
Note that $E(0,\langle i'\rangle)\in Y$ and
\begin{eqnarray*}
&&E(k_{i}\varepsilon_{i}+(k_{j}+1)\varepsilon_{j},\langle
j'\rangle,q)=[E((k_{i}+1)\varepsilon_{i}+(k_{j}+1)\varepsilon_{j},\langle
j'\rangle,q), E(0,i')].
\end{eqnarray*}
We have
 $
E(k_{i}\varepsilon_{i}+(k_{j}+1)\varepsilon_{j},\langle j'\rangle,q)\in
Y$ and
then
\begin{eqnarray*}
&&E(k_{i}\varepsilon_{i}+k_{j}\varepsilon_{j},0,2n+1,q)\\
&\equiv& [E(k_{i}\varepsilon_{i}+(k_{j}+1)\varepsilon_{j},0),
E(0, \langle j'\rangle,2n+1,q)]\\
&=&E(k_{i}\varepsilon_{i}+k_{j}\varepsilon_{j},0,2n+1,q)-
(k_{i}+k_{j}-1)E(k_{i}\varepsilon_{i}+(k_{j}+1)\varepsilon_{j},\langle
j'\rangle,q)\pmod{Y}.
\end{eqnarray*}

 \textit{Subcase 1.1.2}. Suppose $k_{i}=\pi_{i},$ $k_{j}=\pi_{j}.$
By (\ref{e2.3}),
$$E(\pi_{i}\varepsilon_{i}+(\pi_{j}-1)\varepsilon_{j},0,2n+1,q)\in Y.$$
It follows from (\ref{e2.5}) that
$$
E(\pi_{i}\varepsilon_{i}+\pi_{j}\varepsilon_{j},0,2n+1,q)\in
     Y.
$$

\textit{Case 1.2}. Suppose $i\neq j\in \overline{1,n},$  $q=i$ or
$q=j$, say, $q=j.$ Then $k_{j}\neq \pi_{j}.$ It follows from
(\ref{e2.6}) that
$$
E(k_{i}\varepsilon_{i}+k_{j}\varepsilon_{j},0,2n+1,q)\in
     Y.
$$

 Similar to Cases 1.1  and 1.2, using (\ref{e2.7})--(\ref{e2.14}) we may obtain that
$$
E(\varepsilon_{i},\langle j'\rangle,2n+1,q),\quad E(0,\langle i',
j'\rangle,2n+1,q)\in Y \quad\mbox{for all}\; i\neq j\in
\overline{1,n}.
$$
Then by (\ref{e2.15}) and induction one can show that
$S_{3}\subseteq Y.$

\textit{Step 2}. By Part 3 in the proof of Theorem  \ref{t9} and Step 1
we have $S_{4}\setminus \{\pm G\}\subseteq Y.$

\textit{Step 3}. Assert that $\{D_{KO}(f)\mid f\in
\mathfrak{A}_{2}\}\cup S_{2}\subseteq Y.$ Note that
$$\{D_{KO}(f)\mid f\in
\mathfrak{A}_{2}\}\cup S_{2}\setminus \{\pm
E(\pi-\varepsilon_{1},\langle 2',\ldots,n'\rangle,1)\}\subseteq Y$$
can
be easily seen by Lemma \ref{l6}, since $S_{3}\cup S_{4}\setminus
\{\pm G\}\subseteq Y.$
Note that
\begin{eqnarray*}
(3-n)E(\pi-\varepsilon_{1},\langle 2',\ldots,n'\rangle,1)\equiv
[E(\pi-\varepsilon_{1},0,2n+1,1),E(0,\langle
2',\ldots,n'\rangle)]\pmod{Y}
\end{eqnarray*}
and
\begin{eqnarray*}
(4-n)E(\pi-\varepsilon_{1},\langle 2',\ldots,n'\rangle,1)\equiv
[G(\pi-\varepsilon_{1},\langle 2'\rangle,2n+1,1),E(0,\langle
3',\ldots,n'\rangle)]\pmod{Y}.
\end{eqnarray*}
We have $E(\pi-\varepsilon_{1},\langle 2',\ldots,n'\rangle,1)\in Y.$

\textit{Step 4}. Let us show that
\begin{equation}\label{liu9301}
D_{KO}(X(i_{1}, \ldots, i_{r}))\in Y\quad \mbox{for}\; r\notin
\mathfrak{S}_{2}(\lambda, n).
\end{equation}
Note that$$E(0,i_{n}')\in \{D_{KO}(f)\mid f\in
\mathfrak{A}_{2}\}\subseteq Y
$$
and that for $r\in
\overline{0,n},$
$$
E(\pi_{i_{1}}\varepsilon_{i_{1}}+\cdots
+\pi_{i_{r}}\varepsilon_{i_{r}},\langle i_{r+1}', \cdots,
i_{n-1}'\rangle,2n+1,i_{n})\in S_{3}\subseteq Y.
$$
Since  \begin{eqnarray*}\label{e2.20}
&&[E(\pi_{i_{1}}\varepsilon_{i_{1}}+\cdots
 +\pi_{i_{r}}\varepsilon_{i_{r}},\langle i_{r+1}', \cdots,
i_{n-1}'\rangle,2n+1,i_{n}),
E(0,i_{n}')]\nonumber\\
&=&(n\lambda-n+2r+2)D_{KO}(X(i_{1}, \ldots, i_{r})),
\end{eqnarray*}
  (\ref{liu9301}) holds.
\end{proof}

\section{Simplicity and dimension formulas}

Using the spanning set of $\frak{g}$ (Theorem \ref{t10}), let us
prove the following
\begin{theorem}\label{t11}
$\mathfrak{g}$ is a  simple Lie superalgebra.
\end{theorem}
\begin{proof} Let $I$ be a nonzero ideal of $\mathfrak{g}$. By Theorem
\ref{t10}, it suffices to show that $T\cup S\cup
\{D_{KO}(1)\}\subset I.$

\textit{First,} assert that $D_{KO}(1)\in I.$ By Lemma
\ref{l6}, $I$ must contain a nonzero
 element $D_{KO}(f)$ with $\partial_{2n+1}(f)=0.$
 Let
$f= f_{0}x_{n'}+f_{1} $ be the  $x_{n'}$-decomposition. Since
$$
[D_{KO}(f), D_{KO}(x_{n})]=(-1)^{\mathrm{p}(f)}D_{KO}(f_{0}),
$$
 one may assume that $ \partial_{n'}(f)=0.$ Next suppose
$D_{KO}(f):=D_{KO}(f_{0}x_{j'}+f_{1})\in I,$ where
$f=f_{0}x_{j'}+f_{1}$ is the $x_{j'}$-decomposition, $j\in
\overline{1,n-1}$.  Note that
$$
[D_{KO}(f),D_{KO}(x_{j})] =(-1)^{\mathrm{p}(f)}D_{KO}(f_{0})\neq 0.
$$
One may assume that $\partial_{j}(f)=0$ for all $ j\in
\overline{n+1,2n+1}.$ Write
$$
f=a_{0}x^{(k\varepsilon_{i})}+a_{1}x^{((k-1)\varepsilon_{i})}+\cdots
+a_{k-1}x^{(\varepsilon_{i})}+a_{k},
$$
where $i\in \overline{1,n},$ $\partial_{i}(a_{j})=0, j\in
\overline{0,k}.$ Note that
\begin{eqnarray*}\label{eliu7}
&&[D_{KO}(f),D_{KO}(x_{i'})]\nonumber\\
&=&(-1)^{\mu(i)\mathrm{p}(f)}D_{KO}(a_{0}x^{((k-1)\varepsilon_{i})}+a_{1}x^{((k-2)\varepsilon_{i})}+\cdots
+a_{k-1}).
\end{eqnarray*}
The assertion follows.

 \textit{Second},  assert that $T\subseteq I.$ By Lemma \ref{l6},
$E(0,q'),$ $E(\pi_{i},\langle i'\rangle,j)\in I $ for $i,q\in
\overline{1,n}.$ Suppose $i,q\in \overline{1,n}. $
 If $k_{i}\neq \pi_{i}-1$ or $q\neq i$, then
 \begin{equation}\label{liun1}
-E(k_{i}\varepsilon_{i},0,2n+1,q)
=[E(k_{i}\varepsilon_{i}+\varepsilon_{q},0,2n+1,q), E(0,\langle
q'\rangle)]\in I.
   \end{equation}
 If $k_{i}=\pi_{i}-1 $ and $q=i$,  find $j\neq i.$ Then
\begin{eqnarray}\label{liun2}
&&E((\pi_{i}-1)\varepsilon_{i},0,2n+1,i)\nonumber\\
&=&E(\pi_{i},0,2n+1,j)+(n\lambda-\pi_{i}+1)E((\pi_{i}-1)\varepsilon_{i},\langle
i'\rangle,j)\in I.
\end{eqnarray}
The assertion  follows from (\ref{liun1}) and (\ref{liun2}).

\textit{Finally}, it suffices to show that $S\subseteq I.$ This follows directly from  that
$$
-E(0,\langle i'\rangle,2n+1,q)=[E(0,\langle q'\rangle
),E(\varepsilon_{q},\langle i'\rangle,2n+1,q) ]\in I,
$$
where $q\neq i.$
\end{proof}

By Theorems \ref{t5}, \ref{t9} and  \cite [Theorem 4.7]{lh} we can compute the dimension of $\frak{g}$:
\begin{theorem}\label{t12}
\begin{eqnarray*}
\dim
\mathfrak{g}&=&2\Big(\sum_{l=2}^{n}\Big((2^{n-1}-2^{n-l})\sum_{(i_{1},i_{2},\ldots,i_{l})\in
J(l)}
\prod_{c=1}^{l}\pi_{i_{c}}\Big)+\prod_{j=1}^{n}(\pi_{j}+2)\Big)\\
&&-\sum_{k_{i}\in \mathfrak{S}_{2}(\lambda, n)}{n\choose
k_{i}}-2^{n}-\delta'_{n\lambda,-1}.
\end{eqnarray*}
\end{theorem}
\begin{proof}
By \cite [Theorem 4.7]{lh} and (\ref{liue2}), we have
\begin{eqnarray*}
\dim\mathrm{span}_{\mathbb{F}}(S_{1}\cup S_{2})&=&\dim  SHO'(n,n;\underline{t})\\
&=&\sum_{l=2}^{n}\Big((2^{n-1}-2^{n-l})\sum_{(i_{1},i_{2},\ldots,i_{l})\in
J(l)} \prod_{c=1}^{l}\pi_{i_{c}}\Big)\\
&&+\prod_{j=1}^{m}(\pi_{j}+2)-1.
\end{eqnarray*}
By \cite [Theorem 4.7]{lh} and (\ref{e1.1}), we have
\begin{eqnarray*}
\dim \mathrm{span}_{\mathbb{F}}(S_{3}\cup
S_{4})&=&\dim \overline{SHO}(n,n;\underline{t})+1\\
&=&\sum_{l=2}^{n}\Big((2^{n-1}-2^{n-l})\sum_{(i_{1},i_{2},\ldots,i_{l})\in
J(l)}
\prod_{c=1}^{l}\pi_{i_{c}}\Big)\\
&&+\prod_{j=1}^{m}(\pi_{j}+2)-2^{n}.
\end{eqnarray*}
Our formula follows from Theorems \ref{t5}, \ref{t9} and \ref{t8}.
\end{proof}
To  make a comparison between the special odd contact
superalgebras and the other simple  Lie superalgebras of Cartan
type, we list certain known  dimension formulas.

\begin{lemma}[see \cite{lzw,lh}] \label{het3.10}
Suppose  $m,n>2,$ $\underline{t}\in \mathbb{N}^m.$
\begin{enumerate}
\item [$\mathrm{(i)}$]
  $\dim W\left(
m,n;\underline{t}\right) =\left( m+n\right) \cdot 2^n\cdot
p^{\sum_{i=1}^mt_i}.$
\item [$\mathrm{(ii)}$] $ \dim H\left(
m,n;\underline{t}\right) =2^n\cdot p^{\sum_{i=1}^mt_i}-2.$
\item [$\mathrm{(iii)}$] $\dim K\left(
m,n;\underline{t}\right) =\left\{
\begin{array}{ll}
2^n\cdot p^{\sum_{i=1}^mt_i} & \mbox{if}\;\; n-m-3\not\equiv
0 \pmod p \\
 2^n\cdot p^{\sum_{i=1}^mt_i}-1 & \mbox{if}\;\; n-m-3\equiv
0\pmod p.
\end{array}
\right. $
\item [$\mathrm{(iv)}$] $\dim S\left(
m,n;\underline{t}\right) =\left( m+n-1\right) \cdot 2^n\cdot
p^{\sum_{i=1}^mt_i}-m+1.$
\item [$\mathrm{(v)}$] $\dim HO\left(
m,m;\underline{t}\right) =2^m\cdot p^{\sum_{i=1}^mt_i}-1.$
\item [$\mathrm{(vi)}$]
$ \dim SHO(n,n;\underline{t})=
\sum_{l=2}^{n}\Big((2^{n-1}-2^{n-l})\sum\limits_{(i_{1},i_{2},\ldots,i_{l})\in
J(l)}\prod_{c=1}^{l}\pi_{i_{c}}\Big)
+\prod_{j=1}^{n}(\pi_{j}+2)-2^{n}-2. $
\item [$\mathrm{(vii)}$]
$\dim KO(n,n+1;\underline{t})=2^{n+1}p^{\sum_{i=1}^{n}t_{i}}.$
\end{enumerate}
\end{lemma}


\begin{corollary}\label{c1}
$SKO(p+2,p+3;(p-1)/2,\underline{r})$ is not isomorphic to any
Lie superalgebras of Cartan type $W(m,n;\underline{t}),$
$H(m,n;\underline{t}),$ $KO(m,m+1;\underline{t})$ or
$SHO(m,m;\underline{t})$ for arbitrary integers $m,n>2,$
$\underline{r}\in \mathbb{N}^{p+2}$ and $\underline{t}\in
\mathbb{N}^{m}.$
\end{corollary}
\begin{proof}
By Theorem \ref{t12}, we have
\begin{eqnarray*}
\dim
SKO(p+2,p+3;(p-1)/2,\underline{t})&=&2\Big(\sum_{l=2}^{p+2}\Big((2^{p+1}-2^{p+2-l})\sum_{(i_{1},i_{2},\ldots,i_{l})\in
J(l)} \prod_{c=1}^{l}\pi_{i_{c}}\Big)
\\
&&+\prod_{j=1}^{p+2}(\pi_{j}+2)\Big)-2^{p+2}-1.
\end{eqnarray*}
Then $\dim
SKO(p+2,p+3;(p-1)/2,\underline{t})$ is odd and our corollary follows from
  Lemma  \ref{het3.10}.
\end{proof}

\begin{remark} Further comparison between the special odd contact
superalgebras  and the special superalgebras is left to discuss in Section 5,
 where we shall use the structures of the outer superderivation algebras
 (see Corollary \ref{c2} and
Remark \ref{remark5.4}).
\end{remark}
\section{Superderivations}
 Define  a new
multiplication $[\; , \;]$ in $\mathcal{O},$
\begin{eqnarray*}
[a,b]:=D_{KO}(a)(b)-(-1)^{\mathrm{p}(a)}2(\partial_{2n+1}a)b.
\end{eqnarray*}
Then $(\mathcal{O},[\; , \;])$  is  a Lie superalgebra. Since
$\mathfrak{g}''$ is a subalgebra of
$KO(n,n+1;\lambda,\underline{t})$, it is easy to see that
$\{a\mid\mathrm{div}_{\lambda}(a)=0, a\in \mathcal{O}\}$ is a
subalgebra of $\mathcal{O}$.
By \cite{fzj}, the mapping
$$\mathcal{O}\longrightarrow KO(n,n+1;\underline{t}),\quad
D_{KO}: a\longmapsto D_{KO}(a)
$$
 is an isomorphism of Lie superalgebras. Therefore,
$$
\mathfrak{g}''\cong\{a\mid \mathrm{div}_{\lambda}(a)=0, a\in
\mathcal{O}\}.
$$
In this section we  sometimes identify  $D_{KO}(f)$ with $f $ for $f\in
\mathcal {O}.$ Let $\mathrm{Der}\mathfrak{g}$ be the superderivation
algebra of $\mathfrak{g}$. Then
 $\mathrm{Der}\mathfrak{g}$ is a $\mathbb{Z}$-graded Lie superalgebra:
$$
\mathrm{Der}\mathfrak{g}=\oplus_{t\in
\mathbb{Z}}\mathrm{Der}_{t}\mathfrak{g},\;\;\;
\mathrm{Der}_{t}\mathfrak{g}:=\{\phi\in \mathrm{Der}\mathfrak{g}\mid
\phi(\mathfrak{g}_{j})\subseteq \mathfrak{g}_{j+t}, \forall j\in
\mathbb{Z}\}.
$$
For $i\in \overline{1, 2n}$ and $f\in \mathcal {O},$ define
\[ \delta_{if}:=\left\{\begin{array}{ll}0,  &\mbox{
}\partial_{i}(f)=0
\\1,
&\mbox{ }\partial_{i}(f)\neq0. \end{array}\right. \]

Let $T:=\sum_{k=1}^{n}\mathbb{F}h_{k},$ where $h_{k}:=x_{k}x_{k'}.$
\begin{lemma}\label{l12}
Let $h_{i}:=x_{i}x_{i'},$ $i\in \overline{1,n}.$ Then $h_{i}\in
\mathrm{Nor}_{\mathcal{O}}(\mathfrak{g}).$
\end{lemma}
\begin{proof} For $f\in \frak{g},$ we want to show that $[h_{i},f]\in\frak{g}.$ By Theorem \ref{t10},
without loss of generality
one may assume that $f=f_{0}x_{2n+1}+f_{1}$ is an
 element in $T\cup S\cup\{1\}$. Then
\begin{eqnarray}\label{e3.2}
[x_{i}x_{i'},
f]=\delta_{i'f_{0}}f_{0}x_{2n+1}-\delta_{if_{0}}\alpha_{i}f_{0}x_{2n+1}+\delta_{i'f_{1}}f_{1}
-\delta_{if_{1}}(\alpha_{i}+\delta_{i'f_{1}})f_{1},
\end{eqnarray}
where $i\in \overline{1,n}.$ Now it needs only a straightforward verification.
\end{proof}

\begin{lemma}\label{l7}
If $f\in \mathrm{Nor}_{\mathcal {O}}(\mathfrak{\mathfrak{g}})$ and
$\partial_{2n+1}(f)=0.$ Then $f\in \mathfrak{g}''+T.$
\end{lemma}
\begin{proof}
Since $x_{i}\in \mathfrak{g}$ for all $i\in \overline{1,2n},$ we
have $ [x_{i},f]\in \mathfrak{g}$ and then $\Delta([x_{i},f])=0.$ It
is easy to verify that $\partial_{i} (\Delta(f))=0 $ for all $  i\in
\overline{1,2n}. $ It follows that $\Delta (f)\in \mathbb{F}.$
Consequently, $ \Delta(f -\Delta(f)h_{1})=0. $ Therefore, $ f
-\Delta(f)h_{1}\in \mathfrak{g}'' $ and $ f\in \mathfrak{g}''+T. $
\end{proof}
\begin{lemma}\label{p13}
 $ \mathrm{Nor}_{\mathcal{O}}(\mathfrak{g})=\mathfrak{g}''+T.$
\end{lemma}
\begin{proof}
By Lemma \ref{l12}, $\mathfrak{g}''+T\subseteq
\mathrm{Nor}_{\mathcal{O}}(\mathfrak{g}).$  Let us consider the
converse inclusion. Let $f=f_{0}x_{2n+1}+f_{1}\in
\mathrm{Nor}_{\mathcal{O}}(\mathfrak{g})$ be the
$x_{2n+1}$-decomposition. Since $1\in \mathfrak{g},$ we have
$[f_{0}x_{2n+1}+f_{1},1]=2f_{0}\in \mathfrak{g} $ and $f_{0}$ is a
linear combination of $ S_{1}\cup S_{2}\cup \{1\}$. Without loss of
generality, suppose $f_{0}\in S_{1}\cup
S_{2}\cup\{1\}.$\\

\noindent \textit{Case 1}. Suppose $f_{0}\in \mathfrak{A}_{2}\cup
S_{2}$. In this case, $f_{0}$ is integral. Then there is $q\in
\overline{1,n}$ such that
$$f_{0}x_{2n+1}+(-1)^{\mathrm{p}(f_{0})}(n\lambda-\mathrm{zd}(f_{0}))\nabla_{q}(f_{0})\in
\mathfrak{g}''\subseteq \mathrm{Nor}_{\mathcal{O}}(\mathfrak{g}).
$$ So
\begin{eqnarray*}
&&(-1)^{\mathrm{p}(f_{0})}(n\lambda-\mathrm{zd}(f_{0}))\nabla_{q}(f_{0})-f_{1}\\
&=&f_{0}x_{2n+1}+(-1)^{\mathrm{p}(f_{0})}(n\lambda-\mathrm{zd}(f_{0}))\nabla_{q}(f_{0})-(f_{0}x_{2n+1}+f_{1})\in
\mathrm{Nor}_{\mathcal{O}}(\mathfrak{g}). \end{eqnarray*} By Lemma
\ref{l7}, we have
$$
(-1)^{\mathrm{p}(f_{0})}(n\lambda-\mathrm{zd}(f_{0}))\nabla_{q}(f_{0})-f_{1}\in
\mathfrak{g}''+T.
$$
It follows that
\begin{eqnarray*}
f_{0}x_{2n+1}+f_{1}&=&f_{0}x_{2n+1}+f_{1}-(-1)^{\mathrm{p}(f_{0})}(n\lambda-\mathrm{zd}(f_{0}))\nabla_{q}(f_{0})\\
&&+(-1)^{\mathrm{p}(f_{0})}(n\lambda-\mathrm{zd}(f_{0}))\nabla_{q}(f_{0})\\
&=&f_{0}x_{2n+1}+(-1)^{\mathrm{p}(f_{0})}(n\lambda-\mathrm{zd}(f_{0}))\nabla_{q}(f_{0})\\
&&+f_{1}-(-1)^{\mathrm{p}(f_{0})}(n\lambda-\mathrm{zd}(f_{0}))\nabla_{q}(f_{0})\in
\mathfrak{g}''+T.
\end{eqnarray*}

\noindent \textit{Case 2}. Suppose $f_{0}\in \mathfrak{A}_{1}.$  By
(\ref{e13}) there is   $r\in \overline{1,n}$ such that
$f_{0}=X(i_{1},\ldots,i_{n}).$ Since
\begin{eqnarray*}
&&-(n\lambda-n+2r)f_{0}x_{2n+1}+(\mathrm{zd}(f_{1})-2-n\lambda
\delta_{i_{1}f_{1}}+n\lambda\delta_{i_{1}'f_{1}})x_{i_{1}}\partial_{i_{1}}(f_{1})\\&=&[f_{0}x_{2n+1}+f_{1},
x_{2n+1}+n\lambda x_{i_{1}}x_{i_{1}'}] \in \mathfrak{g},
\end{eqnarray*}
we have $n\lambda-n+2r=0,$ that is, $f_{0}x_{2n+1}\in
\frak{g''}\subset \mathrm{Nor}_{\mathcal{O}}(\mathfrak{g}).$
Therefore, $f_{1}\in \mathrm{Nor}_{\mathcal{O}}(g).$ By Lemma
\ref{l7}, we know that $f_{1}\in \mathfrak{g''}+T.$
  Hence
$\mathrm{Nor}_{\mathcal{O}}(\mathfrak{g})\subseteq \mathfrak{g''}+T$
and the proof is complete.
\end{proof}

One can directly verify the following two lemmas.
\begin{lemma}\label{l14}
Let $f\in \mathfrak{g} $ and  $[f,x_{i}]=:b_{i},$ $i\in
\overline{1,2n}.$ Let $\phi\in \mathrm{Der}\mathfrak{g}$ such
that   $\phi(x_{i})=\phi(b_{i})=0$ for all $i\in \overline{1,2n}.$
Then $\phi(f)\in \mathfrak{g}_{-2}.$
\end{lemma}

\begin{lemma}\label{l15}
Let $s\geq -1 $ and  $\phi\in \mathrm{Der}_{t}\mathfrak{g}$
satisfy that $\phi(\mathfrak{g}_{j})=0$ for all $ j\in \overline{-2, s}.$
If $s+t\geq -2,$ then $\phi=0.$
\end{lemma}

\begin{proposition}\label{p16}
$\mathrm{Der}_{-2}\mathfrak{g}=\mathrm{ad}\mathfrak{g}_{-2},$
$\mathrm{Der}_{-3}\mathfrak{g}=0.$
\end{proposition}
\begin{proof} (1) Let us first show that $\mathrm{Der}_{-2}\mathfrak{g}=\mathrm{ad}\mathfrak{g}_{-2}.$
Let $\phi\in \mathrm{Der}_{-2}\mathfrak{g}.$ Then
$\phi(\mathfrak{g}_{j})=0$ for $ j=-2,-1.$ Write
$\phi(x_{2n+1}+n\lambda x_{i}x_{i'})=c_{i}\in \mathfrak{g}_{-2}$ for
$i\in \overline{1,n}.$ Put
$\psi:=\phi+\frac{1}{2}c_{1}\mathrm{ad}1.$ Then
$\psi(\mathfrak{g}_{j})=0$ for $ j=-2,-1 $ and
 $\psi(x_{2n+1}+n\lambda x_{i}x_{i'})=c_{i}-c_{1}$ for all $i\in \overline{1,n}.$ Suppose
$\psi(x_{i}x_{j})=d_{ij}\in \mathfrak{g}_{-2} $ for   $i\neq
j'\in \overline{1,2n}$
 and
    $\psi(x_{i}x_{i'}-x_{j}x_{j'})=e_{ij}\in \mathfrak{g}_{-2}$ for
     $i\neq j\in \overline{1 , n}.$

(1a) Note that
    $$
      [x_{2n+1}+n\lambda x_{i}x_{i'},x_{2n+1}+n\lambda
      x_{1}x_{1'}]=0   \quad \mbox{for all}\; 1\neq i\in \overline{1,n} .
    $$
Applying $\psi$ to the equation above, we have
$$
0=[\psi(x_{2n+1}+n\lambda x_{i}x_{i'}),x_{2n+1}+n\lambda
x_{1}x_{1'}]=-2(c_{i}-c_{1}).
$$
It follows that $c_{i}=c_{1}.$

(1b) Choose  distinct $i, j, k, k'\in \overline{1,2n}.$ Applying
$\psi$ to the following equation
    $$
     [x_{i}x_{j},x_{2n+1}+n\lambda
x_{k}x_{k'}]=0,
    $$
one gets
    $$0=[\psi(x_{i}x_{j}),x_{2n+1}+n\lambda
x_{k}x_{k'}]=[d_{ij},x_{2n+1}+n\lambda x_{k}x_{k'}]=-2d_{ij}.$$ It
follows that $d_{ij}=0.$

(1c) For  distinct $i, j, k\in \overline{1,n},$ applying $\psi$ to
the following equation
    $$
    [x_{i}x_{i'}-x_{j}x_{j'},x_{2n+1}+n\lambda
x_{k}x_{k'}]=0,
    $$
one gets
$$0=[\psi(x_{i}x_{i'}-x_{j}x_{j'}),x_{2n+1}+n\lambda
x_{k}x_{k'}]=-2e_{ij}.$$ It follows that $e_{ij}=0.$

Summarizing, we have $\psi(\mathfrak{g}_{0})=0.$ By Lemma \ref{l15},
$\psi=0,$ that is,
$\mathrm{Der}\mathfrak{g}_{-2}=\mathrm{ad}\mathfrak{g}_{-2}.$

(2) It remains to show that $\mathrm{Der}_{-3}\mathfrak{g}=0.$ Let
$\phi\in \mathrm{Der}\mathfrak{g}_{-3}.$ Then
$\phi(\mathfrak{g}_{j})=0$ for $j=-2,-1,0.$
Write for  $ i\in
\overline{1, 2n},$ $ j\in \overline{1,n}$ with  $j'\neq i,$
$$
\phi(x_{i}x_{2n+1}+(-1)^{\mu(i)}(n\lambda-1)\nabla_{j}(x_{i}))=c_{ij}\in
\mathfrak{g}_{-2}.
$$
Write for $i,j,k\in \overline{1,2n}$ satisfying that $\Delta(x_{i}x_{j}x_{k})=0,$
$$
\phi(x_{i}x_{j}x_{k})=d_{ijk}\in \mathfrak{g}_{-2}.
$$
 Write for $i, k\in \overline{1, n}, j\in \overline{1,2n}$ satisfying that
 $ \Delta(x_{i}x_{i'}x_{j}-x_{k}x_{k'}x_{j})=0,$
$$
\phi(x_{i}x_{i'}x_{j}-x_{k}x_{k'}x_{j})=e_{ijk}\in
\mathfrak{g}_{-2}.
$$

(2a) Fix any  $ i\in \overline{1, 2n} $ and  $ j\in \overline{1, n}$
with $j'\neq i$.
    Applying $\phi$ to the following equation
$$
[x_{i}x_{2n+1}+(-1)^{\mu(i)}(n\lambda-1)\nabla_{j}(x_{i}),x_{2n+1}+n\lambda
x_{j}x_{j'}]=x_{i}x_{2n+1}+(-1)^{\mu(i)}(n\lambda-1)\nabla_{j}(x_{i}),
$$
we have
$$
[\phi(x_{i}x_{2n+1}+(-1)^{\mu(i)}(n\lambda-1)\nabla_{j}(x_{i})),x_{2n+1}+n\lambda
x_{j}x_{j'}]
=\phi(x_{i}x_{2n+1}+(-1)^{\mu(i)}(n\lambda-1)\nabla_{j}(x_{i})).
$$
It follows that $-2c_{ij}=c_{ij},$ that is, $c_{ij}=0.$

(2b) Suppose $\Delta(x_{i}x_{j}x_{k})=0$ for some $i,j,k\in
\overline{1, 2n}$.  Applying $\phi$ to the following equation
$$
[x_{i}x_{j},x_{k}x_{2n+1}+(-1)^{\mu(k)}(n\lambda-1)
x_{i}x_{i'}x_{k}]=\pm(n\lambda-1)x_{i}x_{j}x_{k},
$$
we have
$$
0=\pm(n\lambda-1)\phi(x_{i}x_{j}x_{k})=\pm(n\lambda-1)d_{ijk}.
$$
Thus, if $n\lambda-1 \neq 0$ then $d_{ijk}=0.$ Applying $\phi$ to the following equation
$$
[x_{i}x_{j}x_{k},x_{2n+1}+n\lambda x_{i}x_{i'}]=(\pm
n\lambda+1)x_{i}x_{j}x_{k},
$$
we have
$$
-2d_{ijk}=(\pm n\lambda+1)d_{ijk}.
$$
If $n\lambda-1=0$ then $n\lambda+3\neq 0$. It follows that $d_{ijk}=0.$

(2c) Suppose $\Delta(x_{i}x_{i'}x_{j}-x_{k}x_{k'}x_{j})=0$ for some
$ i, k\in \overline{1,n},$ $ j\in \overline{1,2n}.$ Applying $\phi$
to the following equation
$$
[x_{j}x_{j'}-x_{i}x_{i'},
x_{i}x_{i'}x_{j}-x_{k}x_{k'}x_{j}]=-x_{i}x_{i'}x_{j}+x_{k}x_{k'}x_{j},
$$
we have
$$
0=[x_{j}x_{j'}-x_{i}x_{i'}, e_{ijk}]=-e_{ijk}.
$$
It follows that $e_{ijk}=0.$ The proof is complete.
\end{proof}

\begin{lemma}\label{l17}
Suppose $\phi\in \mathrm{Der}_{-t}\mathfrak{g},$ $ t>3.$ Then
$$\phi(x^{((t-1)\varepsilon_{i})}x_{j})=\phi(x^{((t-1)\varepsilon_{i})}x_{i'}-x^{((t-2)\varepsilon_{i})}x_{k}x_{k'})=0$$
for all $ i, k\in  \overline{1,n} $ and $j\in  \overline{1,2n}$ with $i\neq
j,j', k.$
\end{lemma}
\begin{proof}
Fix any $ i, k\in  \overline{1,n} $ and $j\in  \overline{1,2n}$ with $i\neq
j,j', k.$
Write
$$\phi(x^{((t-1)\varepsilon_{i})}x_{i'}-x^{((t-2)\varepsilon_{i})}x_{k}x_{k'})=d_{ik}\in
\mathfrak{g}_{-2}.
$$
Note that
\begin{eqnarray*}
&&[x^{((t-1)\varepsilon_{i})}x_{i'}-x^{((t-2)\varepsilon_{i})}x_{k}x_{k'},
x_{i}x_{i'}-x_{k}x_{k'}]=
(t-2)(x^{((t-1)\varepsilon_{i})}x_{i'}-x^{((t-2)\varepsilon_{i})}x_{k}x_{k'}),~~~~
\\
&&[x^{((t-1)\varepsilon_{i})}x_{i'}-x^{((t-2)\varepsilon_{i})}x_{k}x_{k'},x_{2n+1}+n\lambda
x_{j}x_{j'}]=(t-2)(x^{((t-1)\varepsilon_{i})}x_{i'}-x^{((t-2)\varepsilon_{i})}x_{k}x_{k'}).
\end{eqnarray*}
Applying $\phi$ to the two equations above, one can obtain that
$(t-2)d_{ik}=0 $ and $td_{ik}=0.$ Since $p>3,$ we have $ d_{ik}=0.$
Given $ i, k\in  \overline{1,n},$ $j\in  \overline{1,2n}$ with $i\neq
j,j', k\neq i,$ apply $\phi$ to the equation below
$$
[x^{((t-1)\varepsilon_{i})}x_{i'}-x^{((t-2)\varepsilon_{i})}x_{k}x_{k'},
x^{(\varepsilon_{i})}x_{j}]=-x^{((t-1)\varepsilon_{i})}x_{j}.
$$
We have $\phi(x^{((t-1)\varepsilon_{i})}x_{j})=0.$
\end{proof}

\begin{lemma}\label{l18}
Suppose $t>3$ and  $\phi\in \mathrm{Der}_{-t}\mathfrak{g}$ such that
 $\phi(x^{(t\varepsilon_{i})})=0$ for all $ i\in
 \overline{1,n}.$ Then $\phi(x^{(k\varepsilon_{i})}x_{j})=0$
 for all $ k\in \overline{0,p^{t_{i}}-1}$ and  $j\in  \overline{1,2n}$ with $i\neq
j'.$
\end{lemma}
\begin{proof}
(1) We first show that $\phi(x^{(k\varepsilon_{i})})=0$ for all $
k\in \overline{0,p^{t_{i}}-1}$. Note that
$\phi(x^{(k\varepsilon_{i})})=0$ for all $ k\in \overline{0,t}.$ If
$k>t,$ by applying $\phi$ to
$[x^{(k\varepsilon_{i})},x_{l}]=\delta_{il'}x^{((k-1)\varepsilon_{i})}$
 for  $l\in \overline{1, 2n} $
 we obtain by inductive hypothesis that
$\phi(x^{((k-1)\varepsilon_{i})})=\phi(x_{l})=0.$ Then by Lemma
\ref{l14}, $\phi(x^{(k\varepsilon_{i})})\in \mathfrak{g}_{-2}\cap
\mathfrak{g}_{k-t-2}.$ Hence $\phi(x^{(k\varepsilon_{i})})=0$ for all
$t<k < p^{t_{i}}.$

(2) Let us show that $\phi(x^{(k\varepsilon_{i})}x_{j})=0$ for all
$k\in \overline{0,p^{t_{i}-1}},$ $ j\in \overline{1,2n}$ and $ i\in
\overline{1,n} $ with $j'\neq i.$

(2.1) Suppose $k<t-1$. As $x^{(k\varepsilon_{i})}x_{j}\in
\mathfrak{g}_{k-1},$ we have $\phi(x^{(k\varepsilon_{i})}x_{j})\in
\mathfrak{g}_{k-1-t}$ and then
$\phi(x^{(k\varepsilon_{i})}x_{j})=0.$

(2.2) By Lemma \ref{l17}, $\phi(x^{((t-1)\varepsilon_{i})}x_{j})=0.$

(2.3) Suppose $k>t-1$ and use induction on $k$. Put
$$
b_{l}:=[x^{(k\varepsilon_{i})}x_{j},x_{l}]=\delta_{il'}x^{((k-1)\varepsilon_{i})}x_{j}
+(-1)^{\mu(i)}\delta_{jl'}x^{(k\varepsilon_{i})}\quad \mbox{for}\;\;l\in
 \overline{1,2n}.
$$
By inductive hypothesis, $\phi(b_{l})=\phi(x_{l})=0.$ Then by Lemma
\ref{l14}, we have $\phi(x^{(k\varepsilon_{i})}x_{j})\in
\mathfrak{g}_{-2}\cap \mathfrak{g}_{k-1-t}$ for all $k>t-1.$
Therefore, $\phi(x^{(k\varepsilon_{i})}x_{j})=0.$ The proof is
complete.
\end{proof}

\begin{lemma}\label{l19}
Suppose $t>3$ and $\phi\in \mathrm{Der}_{-t}\mathfrak{g}$ such that
$\phi(x^{(k\varepsilon_{i})}x_{j})=0$ for all $k\in
\overline{0,p^{t_{i}}-1},$ $i\in  \overline{1,n}$ and  $j\in
\overline{1,2n}$ with $i\neq j'.$ Then
$\phi(x^{(k\varepsilon_{i})}x_{2n+1}+(n\lambda-k)\nabla_{q}(x^{(k\varepsilon_{i})}))=0
$ for all $q\in \widetilde{I}(k\varepsilon_{i},0).$
\end{lemma}
\begin{proof}
(1) Since $t>3$ and $\phi\in \mathrm{Der}_{-t}\mathfrak{g},$ we have
$\phi(x^{(\varepsilon_{i})}x_{2n+1}+(n\lambda-1)\nabla_{q}(x^{(\varepsilon_{i})}))=0.$

(2) Suppose $k=2, t>4$. For any $\phi\in
\mathrm{Der}_{-t}\mathfrak{g},$ it is clear that
$\phi(x^{(2\varepsilon_{i})}x_{2n+1}+(n\lambda-2)\nabla_{q}(x^{(2\varepsilon_{i})}))=0.$

(3) Suppose  $k=2, t=4$. Write
$\phi(x^{(2\varepsilon_{i})}x_{2n+1}+(n\lambda-2)\nabla_{q}(x^{(2\varepsilon_{i})}))=c_{iq}\in
\mathfrak{g}_{-2}.$ Pick $m\neq i,q.$ Applying $\phi$ to the
equation that
    $$
    [x^{(2\varepsilon_{i})}x_{2n+1}+(n\lambda-2)\nabla_{q}(x^{(2\varepsilon_{i})}),
    x_{2n+1}+n\lambda x_{m}x_{m'}]
    =2(x^{(2\varepsilon_{i})}x_{2n+1}+(n\lambda-2)\nabla_{q}(x^{(2\varepsilon_{i})})),
    $$
   one can obtain that $-2c_{iq}=2c_{iq}.$ If follows that 所以$c_{iq}=0.$

(4) Suppose $k>2$ and use induction on $k$. For $j\in
\overline{1,2n},$ put
\begin{eqnarray*}
    b_{j}&:=&[x^{(k\varepsilon_{i})}x_{2n+1}+(n\lambda-k)\nabla_{q}(x^{(k\varepsilon_{i})}),x_{j}]\\
    &=&\delta_{i'j}
    x^{((k-1)\varepsilon_{i})}x_{2n+1}
+\delta_{i'j}(n\lambda-k)\nabla_{q}(x^{((k-1)\varepsilon_{i})})\\
&&+ \delta_{q'j}x^{(k\varepsilon_{i})}x_{q'}+
\delta_{qj}(1-\delta_{ij})(n\lambda-k)x^{(k\varepsilon_{i})}x_{j}+x^{(k\varepsilon_{i})}x_{j}.
\end{eqnarray*}
By induction hypothesis, we have $\phi(x_{j})=\phi(b_{j})=0.$
Then by Lemma \ref{l14},
$$\phi(x^{(k\varepsilon_{i})}x_{2n+1}+(n\lambda-k)\nabla_{q}(x^{(k\varepsilon_{i})}))\in
\mathfrak{g}_{-2}\cap \mathfrak{g}_{k-t}.
$$

(4.1) Suppose $k\neq t-2.$ Then
$\phi(x^{(k\varepsilon_{i})}x_{2n+1}+(n\lambda-k)\nabla_{q}(x^{(k\varepsilon_{i})}))=0.$

(4.2) Suppose $k=t-2.$ Write
$\phi(x^{(k\varepsilon_{i})}x_{2n+1}+(n\lambda-k)\nabla_{q}(x^{(k\varepsilon_{i})}))=c_{iq}\in
\mathfrak{g}_{-2}.$
  Applying $\phi$ to the equation below,
    $$
   [x^{(k\varepsilon_{i})}x_{2n+1}+(n\lambda-k)\nabla_{q}(x^{(k\varepsilon_{i})}),x_{2n+1}+n\lambda x_{q}x_{q'}]
   =k(x^{(k\varepsilon_{i})}x_{2n+1}+(n\lambda-k)\nabla_{q}(x^{(k\varepsilon_{i})}))
    $$
  one gets $-2c_{iq}=kc_{iq}.$

(4.2.1) If $k+2 \not\equiv 0\pmod{p},$ then $c_{iq}=0.$

(4.2.2) If $k+2\equiv 0\pmod{p},$ applying $\phi$ to the equation
\begin{eqnarray*}
&&[x^{((k-1)\varepsilon_{i})}x_{2n+1}+(n\lambda-(k-1))\nabla_{q}(x^{((k-1)\varepsilon_{i})}),
      x^{(\varepsilon_{i})}x_{2n+1}+(n\lambda-1)\nabla_{q}(x^{(\varepsilon_{i})})]\\
      &=&(k-2)(x^{(k\varepsilon_{i})}x_{2n+1}+(n\lambda-k)
     \nabla_{q}(x^{(k\varepsilon_{i})})),
\end{eqnarray*}
we have $c_{iq}=0.$
\end{proof}

\begin{proposition}\label{p20}
If $ t>3$ is not a $p$-power, then
$\mathrm{Der}_{-t}\mathfrak{g}=0.$
\end{proposition}
\begin{proof} Let $
\phi\in \mathrm{Der}_{-t}\mathfrak{g}.$ By Lemmas \ref{l18},
\ref{l19} and Theorem \ref{t10}, it suffices to show that
$\phi(x^{(t\varepsilon_{i})})=0$ for all $i\in \overline{1,n}.$ Let
$\phi(x^{(t\varepsilon_{i})})=c_{i}\in \mathfrak{g}_{-2}$ and pick
$i\neq k\in \overline{1,n}.$ Applying $\phi$ to the following
equation that
$$
[x^{(t\varepsilon_{i})},x_{2n+1}+n\lambda
x_{k}x_{k'}]=(t-2)x^{(t\varepsilon_{i})},
$$
one can obtain that $-2c_{i}=(t-2)c_{i}$ and therefore, $tc_{i}=0.$

If $t\not \equiv0\pmod{p},$ then $c_{i}=0.$ Otherwise, write $t$ as
the $p$-adic from
$$t=\sum_{i=1}^{l}\alpha_{i}p^{i},\;0\leq
\alpha_{i}<
    p,\;\alpha_{l}\neq0.$$
As $t\neq p^{l},$ we have
   \begin{eqnarray*}
&&\phi(x^{(p^{l}\varepsilon_{i})}x_{i'}-x^{((p^{l}-1)\varepsilon_{i})}x_{k}x_{k'})\in
\mathfrak{g}_{[p^{l}-1-t]}=0,
   \\
   &&
\phi(x^{((t-p^{l}+1)\varepsilon_{i})})\in
\mathfrak{g}_{[-p^{l}-1]}=0.
   \end{eqnarray*}
 Applying $\phi$ to
    $$
[x^{((t-p^{l}+1)\varepsilon_{i})},
x^{(p^{l}\varepsilon_{i})}x_{i'}-x^{((p^{l}-1)\varepsilon_{i})}x_{k}x_{k'}]={t\choose
p^{l}}x^{(t\varepsilon_{i})},
    $$
  one gets $\phi(x^{(t\varepsilon_{i})})=0.$
The proof is complete.
\end{proof}

\begin{proposition}\label{p21}
Let $t= p^{d} $ for some $d\in \mathbb{N}.$
Then
$\mathrm{Der}_{-t}\mathfrak{g}=\mathrm{span}_{\mathbb{F}}\{\partial_{i}^{t}\mid
i\in \overline{1,n}\}.$
\end{proposition}
\begin{proof} Let $\phi\in \mathrm{Der}_{-t}\mathfrak{g}.$
Then $\phi(x^{(t\varepsilon_{i})})=c_{i}\in \mathfrak{g}_{-2}.$
Putting $\psi:=\phi-\sum_{j=1}^{n}c_{j}\partial_{j}^{t},$ we have
$\psi(x^{(t\varepsilon_{i})})=0.$ Then by Lemmas \ref{l18},
\ref{l19} and Theorem \ref{t10},  we have $\psi=0,$ that is,
$\phi=\sum_{j=1}^{n}c_{j}\partial_{j}^{t}. $ This proves
$$\mathrm{Der}\mathfrak{g}_{-t}\subseteq\mathrm{span}_{\mathbb{F}}\{\partial_{i}^{t}\mid i\in \overline{1,n}\}.$$
Note that in general,
$$[\partial_{i}^{p^{r}}, D_{KO}(a)]=D_{KO}(\partial_{i}^{p^{r}}(a))\quad \mbox{and}
\quad \mathrm{div}_{\lambda}(\partial_{i}^{p^{r}}(a))=\partial_{i}^{p^{r}}(\mathrm{div}_{\lambda}(a))
$$ for $ a\in\mathcal{O}$ and $ r\in \mathbb{N}.$ One sees that
 $\partial_{i}^{p^{r}}$ is a derivation of $\frak{g}$ and then the
 converse inclusion holds.
\end{proof}

Given $i\in \overline{1,n},$ an element $a\in \mathcal{O} $ is
called $x_{i}$-truncated if each nonzero monomial
$kx^{(\alpha)}x^{u}$ of $a$ satisfies  that $\alpha_{i}<
p^{t_{i}}-1$. We need a technical lemma:

\begin{lemma}\label{l22}
Suppose $a_{i}\in \mathcal{O}$ satisfy that
 $\partial_{i}(a_{j})=(-1)^{\mu(i)\mu(j)}\partial_{j}(a_{i}) $ for all
$ i,j\in  \overline{1,2n+1} .$ The the following statements hold.

   $(1)$  $\partial_{i}(a_{i})=0 $ for all $i\in  \overline{n+1,2n+1};$

    $(2)$   For $i\in \overline{1,n},$
    $a_{i}$ has only one nonzero monomial of the form
     $c_{i}x^{(\pi_{i}\varepsilon_{i})}, c_{i}\in \mathbb{F};$

    $(3)$   If $a_{i}$ is $x_{i}$-truncated, then there exists $f\in
    \mathcal{O}$ such that $a_{i}=\partial_{i}(f)$
     for all $  i\in \overline{1,2n+1}.$
\end{lemma}
\begin{lemma}\label{l23}
If $\phi\in \mathrm{Der}\mathfrak{g},$ then there is
$f\in
    \mathcal{O}$
   such that $(\phi-\mathrm{ad}f)(\mathfrak{g}_{j})=0$ for  $j=-1,-2.$
\end{lemma}
\begin{proof} Let $\phi\in \mathrm{Der}_{\alpha}\mathfrak{g},$
$\alpha\in \mathbb{Z}_{2}.$ Put
$$
a_{2n+1}:=(-1)^{\alpha}\frac{1}{2}\phi(1)
$$
and for  $i\in \overline{1, 2n},$
$$
a_{i}:=(-1)^{\mu(i)\alpha+\mu(i)}\phi(x_{i'})+(-1)^{\mu(i')}x_{i'}a_{2n+1}.
$$
By applying $\phi$ to $[1,x_{i'}]=0$ for  $i\in \overline{1, 2n},$
we have

$$
(-1)^{\alpha}2[a_{2n+1},
x_{i'}]+(-1)^{\alpha+\mu(i)\alpha+\mu(i)}[1,
a_{i}]+(-1)^{\alpha+\mu(i)\alpha}[1, a_{2n+1}x_{i'}]=0,
$$
that is,
\begin{eqnarray*}
&& (-1)^{\alpha+\mu(i)\alpha}2\partial_{i}(a_{2n+1})+2\partial_{2n+1}(a_{2n+1}x_{i'})\\
&=&
 (-1)^{\alpha+\mu(i)\alpha+\mu(i)}
2\partial_{2n+1}(a_{i})
+(-1)^{\alpha+\mu(i)\alpha}2\partial_{2n+1}(a_{2n+1}x_{i'}).
\end{eqnarray*}
It is easy to see that $\partial_{2n+1}(a_{2n+1})=0.$
Therefore,
$$
\partial_{i}(a_{2n+1})=(-1)^{\mu(i)}\partial_{2n+1}(a_{i})\quad
\mbox{for all} \; i\in \overline{1,2n}.
$$
Since $[x_{i'},x_{j'}]=0$  for $i,j\in \overline{1,2n}$ with $j\neq
i',$ we have
$$
[\phi(x_{i'}), x_{j'}]+(-1)^{\alpha\mu(i)}[x_{i'}, \phi(x_{j'})]=0.
$$
It follows that
\begin{eqnarray*}
 &&(-1)^{\mu(i)\alpha+\mu(i)+\mu(j)(\alpha+\mu(i'))}\partial_{j}(a_{i})
-(-1)^{\mu(i)\alpha+\mu(i)+\alpha+\mu(i')}\partial_{2n+1}(a_{i})x_{j'}\\
&&-(-1)^{\mu(i)\alpha+\mu(j)(\alpha+\mu(i'))}\partial_{j}(x_{i'}a_{2n+1})
+(-1)^{\mu(i)\alpha+\alpha+\mu(i')}\partial_{2n+1}(x_{i'}a_{2n+1})x_{j'}\\
&&
+(-1)^{\alpha(\mu(i)+\mu(j))+\mu(j)+\mu(i')}\partial_{i}(a_{j})-
(-1)^{\alpha(\mu(i)+\mu(j))+\mu(j)}x_{i'}\partial_{2n+1}(a_{j})\\
&&
+(-1)^{\alpha(\mu(i)+\mu(j))+\mu(i')}\partial_{i}(x_{j'}a_{2n+1})-
(-1)^{\alpha(\mu(i)+\mu(j))}x_{i'}\partial_{2n+1}(x_{j'}a_{2n+1})=0.
\end{eqnarray*}
Note that
\begin{eqnarray*}
&&-(-1)^{\mu(i)\alpha+\mu(i)+\alpha+\mu(i')}\partial_{2n+1}(a_{i})x_{j'}
+(-1)^{\alpha(\mu(i)+\mu(j))+\mu(i')}\partial_{i}(x_{j'}a_{2n+1})=0;
\\
&&
-(-1)^{\mu(i)\alpha+\mu(j)(\alpha+\mu(i'))}\partial_{j}(x_{i'}a_{2n+1})-
(-1)^{\alpha(\mu(i)+\mu(j))+\mu(j)}x_{i'}\partial_{2n+1}(a_{j})=0;
\\
&&
(-1)^{\mu(i)\alpha+\alpha+\mu(i')}\partial_{2n+1}(x_{i'}a_{2n+1})x_{j'}=0;
\\
&&
(-1)^{\alpha(\mu(i)+\mu(j))}x_{i'}\partial_{2n+1}(x_{j'}a_{2n+1})=0.
\end{eqnarray*}
We have
$$
\partial_{j}(a_{i})=(-1)^{\mu(i)\mu(j)}\partial_{i}(a_{j}).
$$
In the same way, applying $\phi$ to $[x_{i},
x_{i'}]=(-1)^{\mu(i)}1$  for  $i\in \overline{1,2n},$ we
conclude that
\begin{eqnarray*}
&&(-1)^{\mu(i')\alpha+\mu(i')+\mu(i)(\alpha+\mu(i))}\partial_{i}(a_{i'})-
(-1)^{\mu(i')\alpha+\mu(i')+\alpha+\mu(i)}\partial_{2n+1}(a_{i'})x_{i'}
\\
&&+(-1)^{\mu(i')\alpha+\mu(i)(\alpha+\mu(i))}\partial_{i}(x_{i}a_{2n+1})-
(-1)^{\mu(i')\alpha+\alpha+\mu(i)}\partial_{2n+1}(x_{i}a_{2n+1})x_{i'}
\\
&&+(-1)^{\alpha\mu(i')+\mu(i)\alpha+\mu(i)+\mu(i)}\partial_{i'}(a_{i})-
(-1)^{\alpha\mu(i')+\mu(i)\alpha+\mu(i)}x_{i}\partial_{2n+1}(a_{i})
\\
&&+(-1)^{\alpha\mu(i')+\mu(i)\alpha+\mu(i)}\partial_{i'}(x_{i'}a_{2n+1})-
(-1)^{\alpha\mu(i')+\mu(i)\alpha}x_{i}\partial_{2n+1}(x_{i'}a_{2n+1})
\\
&=&(-1)^{\mu(i)+\alpha}2a_{2n+1}.
\end{eqnarray*}
Since
\begin{eqnarray*}
&&-(-1)^{\mu(i')\alpha+\mu(i')+\alpha+\mu(i)}\partial_{2n+1}(a_{i'})x_{i'}
+(-1)^{\alpha\mu(i')+\mu(i)\alpha+\mu(i)}\partial_{i'}(x_{i'}a_{2n+1})
\\
&&
+(-1)^{\mu(i')\alpha+\mu(i)(\alpha+\mu(i))}\partial_{i}(x_{i}a_{2n+1})-
(-1)^{\alpha\mu(i')+\mu(i)\alpha+\mu(i)}x_{i}\partial_{2n+1}(a_{i})
\\
&=&(-1)^{\mu(i)+\alpha}2a_{2n+1}
\end{eqnarray*}
and
$$
(-1)^{\mu(i')\alpha+\alpha+\mu(i)}\partial_{2n+1}(x_{i}a_{2n+1})x_{i'}=(-1)^{\alpha\mu(i')+\mu(i)\alpha}x_{i}\partial_{2n+1}(x_{i'}a_{2n+1})=0,
$$
one can gets
$$
\partial_{i}(a_{i'})=(-1)^{\mu(i)\mu(i')}\partial_{i'}(a_{i})\quad \mbox{for all} \; i\in \overline{1,2n}.
$$
We have proved that
$$\partial_{i}(a_{j})=(-1)^{\mu(i)\mu(j)}\partial_{j}(a_{i}) \quad \mbox{for
all}\; i,j\in \overline{1,2n+1}.$$

Next, we want to show that $a_{i}$ is $x_{i}$-truncated for  $i\in
\overline{1, n}.$ Pick $i,i'\neq j\in \overline{1,n}.$
    Applying $\phi$ to the following equation that
    $$
[x_{2n+1}+n\lambda x_{j}x_{j'},1]=2,
    $$
one gets
$$
[\phi(x_{2n+1}+n\lambda x_{j}x_{j'}),1]+[x_{2n+1}+n\lambda
x_{j}x_{j'},\phi(1)]=2\phi(1),
$$
that is,
\begin{eqnarray*}
&&-(-1)^{\alpha}2\partial_{2n+1}(\phi(x_{2n+1}+n\lambda
x_{j}x_{j'}))-\mathfrak{D}(\phi(1))\\
&&-2x_{2n+1}\partial_{2n+1}(\phi(1))+n\lambda(x_{j'}\partial_{j'}(\phi(1))
-x_{j}\partial_{j}(\phi(1)))+2\phi(1)\\
&=&2\phi(1).
\end{eqnarray*}
Hence
\begin{eqnarray}\label{e3.3}
-\partial_{2n+1}(\phi(x_{2n+1}+n\lambda
x_{j}x_{j'}))=\mathfrak{D}(a_{2n+1})+n\lambda(x_{j}\partial_{j}(a_{2n+1})-x_{j'}\partial_{j'}(a_{2n+1})).
\end{eqnarray}
Applying $\phi$ to the equation that
$$
[x_{i'},x_{2n+1}+n\lambda x_{j}x_{j'}]=-x_{i'},
$$
one can get
$$
[\phi(x_{i'}),x_{2n+1}+n\lambda
x_{j}x_{j'}]+[x_{i'},\phi(x_{2n+1}+n\lambda
x_{j}x_{j'})]=-\phi(x_{i'}),
$$
that is,
\begin{eqnarray*}
&&[(-1)^{\mu(i)\alpha+\mu(i)}a_{i}-(-1)^{\mu(i)\alpha}x_{i'}a_{2n+1},x_{2n+1}+n\lambda
x_{j}x_{j'}]\\
&=&-[x_{i'},\phi(x_{2n+1}+n\lambda
x_{j}x_{j'})]-(-1)^{\mu(i)\alpha+\mu(i)}a_{i}+(-1)^{\mu(i)\alpha}x_{i'}a_{2n+1}.
\end{eqnarray*}
In combination with (\ref{e3.3}), we have
\begin{eqnarray*}
&&(-1)^{\mu(i)\alpha+\mu(i)}\mathfrak{D}(a_{i})-2(-1)^{\mu(i)\alpha+\mu(i)}a_{i}+2(-1)^{\mu(i)\alpha+\mu(i)}\partial_{2n+1}(a_{i})x_{2n+1}\\
&&+(-1)^{\mu(i)\alpha+\mu(i)}n\lambda(\partial_{j}(a_{i})x_{j}-\partial_{j'}(a_{i})x_{j'})
-(-1)^{\mu(i)\alpha}\mathfrak{D}(x_{i'}a_{2n+1})\\
&&+(-1)^{\mu(i)\alpha}2x_{i'}a_{2n+1}+n\lambda(-(-1)^{\mu(i)\alpha}\partial_{j}(x_{i'}a_{2n+1}))x_{j}
-(-1)^{\mu(i)\alpha+\mu(i)}\partial_{j'}(a_{i})x_{j'})\\
&&-\partial_{i}(\phi(x_{2n+1}+n\lambda
x_{j}x_{j'}))+x_{i'}(\mathfrak{D}(a_{2n+1})+n\lambda(x_{j}\partial_{j}(a_{2n+1})-x_{j'}\partial_{j'}(a_{2n+1})))\\
&=&-(-1)^{\mu(i)\alpha+\mu(i)}a_{i}+(-1)^{\mu(i)\alpha}x_{i'}a_{2n+1}.
\end{eqnarray*}
By Lemma \ref{l22}(2), $a_{i}=x+c_{i}x^{(\pi_{i}\varepsilon_{i})},$
where $x$ is $x_{i}$-truncated. Substituting  $a_{i}$ with
$x+c_{i}x^{(\pi_{i}\varepsilon_{i})}$ and observing   the
coefficient of $x^{(\pi_{i}\varepsilon_{i})}$ we obtain that
$\pi_{i}c_{i}-2c_{i}=-c_{i}.$ Consequently, $2c_{i}=0,$ that is,
$a_{i}$ is $x_{i}$-truncated. By Lemma \ref{l22}(1) and  (3), there
exists $f\in \mathcal {O}$ such that
$$
a_{i}=\partial_{i}(f) \quad \mbox{for all}\; i\in
\overline{1,2n+1}.
$$
Then we have
$$
(\phi-\mathrm{ad}(f))(\mathfrak{g}_{j})=0\quad \mbox{for}\;
j=-2,-1.
$$
By Lemma \ref{l15}, the proof is complete.
\end{proof}

\begin{lemma}\label{p24}
Suppose $\phi\in \mathrm{Der}_{t}\mathfrak{g},$ $t\geq-1.$ Then
there exists $f\in \mathcal{O}$ such that $\phi=\mathrm{ad}f.$
\end{lemma}
\begin{proof} Let $\phi\in \mathrm{Der}\mathfrak{g}.$
By Lemma \ref{l23}, there exists $f\in \mathcal{O}$ such that
$(\phi-\mathrm{ad}(f))(\mathfrak{g}_{j})=0 $ for $j=-2,-1.$ Since $t\geq-1,$
by Lemma \ref{l15} we have $\phi-\mathrm{ad}f=0,$ that is,
$\phi=\mathrm{ad}f.$
\end{proof}

We are in the position to determine completely  the superderivation algebra of
 the special odd contact superalgebra $\frak{g}.$
 Recall that $T:=\sum_{k=1}^{n}\mathbb{F}x_{k}x_{k'}.$
\begin{theorem}\label{t25}
$\mathrm{Der}\mathfrak{g}=\mathrm{ad}\mathfrak{g}''+\mathrm{ad}T\oplus
\mathrm{span}_{\mathbb{F}}\{\partial_{i}^{p^{d_{i}}}\mid 1\leq
d_{i}\leq t_{i}-1, i\in \overline{1,n}\}.$
\end{theorem}
\begin{proof}
Note that any $p$-power of an even superderivation is  again   an
even superderivation and that
$[[\frak{g}'',\frak{g}''],[\frak{g}'',\frak{g}'']]=\frak{g}.$ By
Lemma \ref{p13}, we have
$\mathrm{Der}\mathfrak{g}\supseteq\mathrm{ad}\mathfrak{g}''+\mathrm{ad}T\oplus
\mathrm{span}_{\mathbb{F}}\{\partial_{i}^{p^{d_{i}}}\mid 1\leq
d_{i}\leq t_{i}-1, i\in \overline{1,n}\}.$ The converse inclusion
follows from  Lemma \ref{p13}, Propositions \ref{p16}, \ref{p20}, \ref{p21}
and Lemma \ref{p24}.
\end{proof}

\begin{remark} The degree derivation $\deg$ of $\frak{g} $ is
$$-\mathrm{ad}x_{2n+1}=-\mathrm{ad}(x_{2n+1}+n\lambda x_{1}x_{1'})+n\lambda \mathrm{ad}(x_{1}x_{1'})\in
 \mathrm{ad}\mathfrak{g}''+\mathrm{ad}T.$$

\end{remark}
\section{First cohomology}
Recall that  $\mathfrak{g}$ denotes the special odd contact
superalgebra $SKO(n,n+1;\lambda,\underline{t})$, where $n\geq 3$ is
an integer and  $\underline{t}:=(t_1,\ldots,t_n)$ is an $n$-tuple of
positive integers. Denote the outer superderivation algebra of
$\mathfrak{g}$ by $ \mathrm{Der}_{\mathrm{out}}\mathfrak{g}:=
\mathrm{Der}\mathfrak{g}/\mathrm{ad}\mathfrak{g},$ which is isomorphic to the first cohomology group.  In this section
we shall characterize the structure of
$\mathrm{Der}_{\mathrm{out}}\mathfrak{g} $ and give an application for the isomorphism problem between modular Lie
superalgebras of Cartan type.

For simplicity, we introduce some symbols. Let $(i_1,i_2,\ldots,i_k)$ be a $k$-tuple of arbitrary positive
integers. As in usual, put
$$\mathrm{sgn}(i_1,i_2,\ldots,i_k):=\frac{\prod_{1\leq{j}<l\leq k}(i_{l}-i_{j} )}
{|\prod_{1\leq{j}<l\leq k}(i_{l}-i_{j} )|},$$ which is $\pm1.$ Write
the  integers
\begin{eqnarray*}\label{liuy51}
&&l_{\bar{0}}(\lambda,n):=\sum_{k\in \mathfrak{S}_{0}(\lambda,
n)\atop n-k \;\mathrm{is}\; \mathrm{even}}{n\choose k}+\sum_{k\in \mathfrak{S}_{2}(\lambda,
n)\atop n-k\; \mathrm{is}\; \mathrm{odd}}{n\choose k},\\
&&l_{\bar{1}}(\lambda,n):=\sum_{k\in \mathfrak{S}_{0}(\lambda,
n)\atop n-k\; \mathrm{is}\; \mathrm{odd}}{n\choose k}+\sum_{k\in \mathfrak{S}_{2}(\lambda,
n)\atop n-k\; \mathrm{is}\; \mathrm{even}}{n\choose k}.
\end{eqnarray*}
Here we recall that for  $\lambda\in \mathbb{F}$,
$$
\mathfrak{S}_{l}(\lambda, n):=\{k\in \overline{0,n}\mid
n\lambda-n+2k+l=0 \in\mathbb{F}\}.
$$

Let $V:=V_{\bar{0}}\oplus V_{\bar{1}}$ be a $\mathbb{Z}_{2}$-graded
vector space
 where
\begin{eqnarray*}
&&
V_{\bar{0}}:=V_{00}\oplus V_{01}\oplus V_{02}\oplus V_{03},\;
V_{03}:=\delta'_{n\lambda,-1}\mathbb{F}\cdot 1,\; V_{\bar{1}}:=V_{11}\oplus V_{12};\\
&&V_{01}:=\mathrm{span}_{\mathbb{F}}\{X_{i_1,\ldots,i_r}\mid r\in\mathfrak{S}_{2}(\lambda,
n),(i_{1}, \ldots, i_{r})\in \mathbf{J}(r), n-r\;\mbox{ is odd}\};
\\
&&V_{02}:=\mathrm{span}_{\mathbb{F}}\{Y_{j_1,\ldots,j_l}\mid l\in\mathfrak{S}_{0}(\lambda,
n),(j_{1}, \ldots, j_{l})\in \mathbf{J}(l), n-l\;\mbox{ is even}\};\\
&&V_{11}:=\mathrm{span}_{\mathbb{F}}\{X_{i_1,\ldots,i_r}\mid r\in\mathfrak{S}_{2}(\lambda,
n),(i_{1}, \ldots, i_{r})\in \mathbf{J}(r), n-r\;\mbox{ is even}\};
\\
&&V_{12}:=\mathrm{span}_{\mathbb{F}}\{Y_{j_1,\ldots,j_l}\mid l\in\mathfrak{S}_{0}(\lambda,
n),(j_{1}, \ldots, j_{l})\in \mathbf{J}(l), n-l\;\mbox{ is odd}\};\\
&&\dim V_{00}=|t|-n,\; \dim V_{01}\oplus V_{02}=l_{\bar{0}}(n,\lambda), \;
\dim V_{\bar{1}}=l_{\bar{1}}(n,\lambda).
\end{eqnarray*}
 Note that $V$ is of dimension
$$l_{\bar{0}}(\lambda,n)+l_{\bar{1}}(\lambda,n)+|t|-n+\delta'_{n\lambda,-1}.$$
Moreover,$V$ becomes a
  Lie superalgebra
by letting
$$[V_{00},V]=
[V_{03},V]=[V_{01}+V_{11},V_{01}+V_{11}]=[V_{02}+V_{12},V_{02}+V_{12}]=0
$$
and
\begin{eqnarray*}
[X_{i_1,\ldots,i_r},Y_{j_1,\ldots,j_l}]&=&(-1)^{r(n-r+1)}[Y_{j_1,\ldots,j_l},X_{i_1,\ldots,i_r}]\\
&=&\left\{\begin{array}{ll}0,  &\mbox{ }(i_1,\ldots,i_r)\neq
(j_{l+1},\ldots,j_n)
\\\delta'_{n\lambda,-1}\mathrm{sgn}(i'_{r+1},\ldots,i'_n,i'_2,\ldots,i'_r)1,
&\mbox{ }(i_1,\ldots,i_r)= (j_{l+1},\ldots,j_n).
\end{array}\right.
\end{eqnarray*}
In particular, if $\delta'_{n\lambda,-1}=0 $ then $V$ is an abelian Lie
algebra.

 Let $f\in \mathrm{End} V$ be such that
$f|_{V_{00}}=0,$ $f|_{V_{01}\oplus V_{02}\oplus
V_{\bar{1}}}=\mathrm{id}_{V_{01}\oplus V_{02}\oplus V_{\bar{1}}}$ and
$f|_{V_{03}}=2\mathrm{id}_{V_{03}}.$
 Denote the semidirect product by  $\mathfrak{L}:=\mathbb{F}f\ltimes V.$ Then
$\mathfrak{L}$ is a  Lie superalgebra of dimension
$l_{\bar{0}}(\lambda,n)+l_{\bar{1}}(\lambda,n)+|t|-n+1+\delta'_{n\lambda,-1}$.
Consider the centralizer of $\mathfrak{g}$ in $KO$,
$$
C_{KO}(\mathfrak{g}):=\big\{e\in KO\mid [e,\mathfrak{g}]=0\big\}.
$$
\begin{lemma}\label{l26}
$C_{KO}(\mathfrak{g})=0.$
\end{lemma}
\begin{proof}
Let $f$ be an arbitrary element of $C_{KO}(\mathfrak{g})$ and $f=
f_{0}x_{2n+1}+f_{1} $ be the $x_{2n+1}$-decomposition. Note that
$\frak{g}$ contains $1,$ $x_{i}$ for $ i\in \overline{1,2n}$ and
$x_{2n+1}+n\lambda x_{1}x_{1'}.$ We have $f_{0}=0$ since
$2f_{0}=[f,1]=0.$ Similarly, one gets $f\in \mathbb{F}$ from the
equation that $[f,x_{i}]=0$ for  $i\in \overline{1,2n}.$ Since
$-2f=[f,x_{2n+1}+n\lambda x_{1}x_{1'}]=0,$ we have $f=0.$
\end{proof}
\begin{theorem}\label{t26}
$\mathrm{Der}_{\mathrm{out}}\mathfrak{g}\cong \mathfrak{L} $  and
$$\dim\mathrm{Der}_{\mathrm{out}}\mathfrak{g}=\sum_{k\in \mathfrak{S}_{0}(\lambda, n)}{n\choose
k}+\sum_{k\in \mathfrak{S}_{2}(\lambda, n)}{n\choose
k}+|t|-n+1+\delta'_{n\lambda,-1}.
$$
\end{theorem}
\begin{proof}
Let $\rho: \mathrm{Der}\mathfrak{g}\longrightarrow
\mathrm{Der}_{\mathrm{out}}\mathfrak{g}:=\mathrm{Der}\mathfrak{g}/\mathrm{ad}\mathfrak{g}$
be the canonical homomorphism. By Theorem \ref{t25},
$\mathrm{Der}_{\mathrm{out}}\mathfrak{g}$ is spanned by
\begin{eqnarray}\label{e27}
&&\{\rho(\mathrm{ad}X(i_1,\ldots,i_r))\mid
r\in\mathfrak{S}_{2}(\lambda, n),(i_{1}, \ldots, i_{r})\in
\mathbf{J}(r)\}\cup\{\rho(\mathrm{ad}y)\mid y\in
S_{5}\}\cup\{\rho(\mathrm{ad}x_{1}x_{1'})\}\nonumber\\
&&\cup \{\rho(\partial_{i}^{p^{d_{i}}})\mid 1\leq d_{i}\leq t_{i}-1,
i\in
\overline{1,n}\}\cup\{\delta'_{n\lambda,-1}\rho(\mathrm{ad}G)\}.
\end{eqnarray}
Note that
$\mathrm{span}_{\mathbb{F}}\{\rho(\partial_{i}^{p^{d_{i}}})\mid
1\leq d_{i}\leq t_{i}-1, i\in \overline{1,n}\}\cap
\rho(\mathrm{ad}\mathfrak{g}'')=0$ and
$\mathbb{F}(\rho(\mathrm{ad}x_{1}x_{1'}))\cap\mathfrak{\rho(\mathrm{ad}g'')}=0.$
We may prove from Lemma \ref{l26} that the  above  set (\ref{e27}) is
$\mathbb{F}$-linearly independent. Suppose
$$a\in
\{X(i_1,\ldots,i_r)\mid r\in\mathfrak{S}_{2}(\lambda, n),(i_{1},
\ldots, i_{r})\in \mathbf{J}(r)\}\cup S_{5}.
$$ For $i,j\in
\overline{1,n},$ $1\leq d_{i}\leq t_{i}-1$, $1\leq d_{j}\leq
t_{j}-1,$ $(i_{1}, \ldots, i_{r})\in \mathbf{J}(r)$ and $(j_{1},
\ldots, j_{l})\in \mathbf{J}(l),$ we have
\begin{eqnarray*}
&& [\rho(\mathrm{ad}x_{1}x_{1'}),
\rho(\mathrm{ad}a)]=\rho(\mathrm{ad}a),\\
&&[\rho(\mathrm{ad}x_{1}x_{1'}),\rho(\mathrm{ad}G)]=2\rho(\mathrm{ad}G),\\
&&[\rho(\partial_{i}^{p^{d_{i}}}),\rho(\mathrm{Der}(\mathfrak{g}))]=[\rho(\mathrm{ad}S_{5}),\rho(\mathrm{ad}S_{5})]=0,\\
&&[\rho(\mathrm{ad}a),
\rho(\mathrm{ad}G)]=[\rho(\mathrm{ad}X(i_{1},\ldots,i_{r}),\rho(\mathrm{ad}X(j_{i},\ldots,j_{l}))]=0,\\
&&[\rho(\mathrm{ad}X(i_1,\ldots,i_r)),\rho(\mathrm{ad}E(\pi_{j_1}+\cdots+\pi_{j_l},\langle
j'_{l+1},\ldots,j'_n\rangle,2n+1))]\\
&=&(-1)^{r(n-r+1)}[\rho(\mathrm{ad}E(\pi_{j_1}+\cdots+\pi_{j_l},\langle j'_{l+1},\ldots,j'_n\rangle,2n+1)),\rho(\mathrm{ad}X(i_1,\ldots,i_r))]\\
&=&\left\{\begin{array}{ll}0,  &\mbox{ }(i_1,\ldots,i_r)\neq
(j_{l+1},\ldots,j_n)
\\\delta'_{n\lambda,-1}\mathrm{sgn}(i'_{r+1},\ldots,i'_n,i'_2,\ldots,i'_r)\rho(\mathrm{ad}G), &\mbox{ }(i_1,\ldots,i_r)=
(j_{l+1},\ldots,j_n).
\end{array}\right.
\end{eqnarray*}
Now one may easily establish a  Lie superalgebra isomorphism of
$\mathrm{Der}_{\mathrm{out}}\mathfrak{g}$ onto $\mathfrak{L} $.  The
dimension formula follows.
\end{proof}
\begin{corollary}\label{c2}
$SKO(p+2,p+3;(p-1)/2,\underline{r})$ is not isomorphic to
$S(m,n;\underline{t}),$ $HO(m,m;\underline{t})$ or
$K(m,n;\underline{t})$ for any $m,n>2$, $\underline{r}\in
\mathbb{N}^{p+2}$ and $\underline{t}\in \mathbb{N}^{m}.$
\end{corollary}
\begin{proof} By Theorem \ref{t26}, \cite [Theorems
2.12 and 2.4(ii)]{lz3} and \cite [Theorem 22]{lzw}, one sees that
\begin{eqnarray*}
&&
\mathrm{Der}_{\mathrm{out}}S(m,n;\underline{t})\not\cong\mathrm{Der}_{\mathrm{out}}SKO(p+2,p+3;
(p-1)/{2},\underline{r}),\\
&&
\mathrm{Der}_{\mathrm{out}}HO(m,m;\underline{t})\not\cong\mathrm{Der}_{\mathrm{out}}
SKO(p+2,p+3;(p-1)/{2},\underline{r}),
\\&&
\mathrm{Der}_{\mathrm{out}}K(m,n;\underline{t})\not\cong\mathrm{Der}_{\mathrm{out}}SKO(p+2,p+3;(p-1)/2,
\underline{r})
\end{eqnarray*}
for all $m,n>2$, $\underline{r}\in \mathbb{N}^{p+2}$ and
$\underline{t}\in \mathbb{N}^{m}.$ Therefore, $SKO(p+2,p+3;
(p-1)/2,\underline{r})$ is not isomorphic to $S(m,n;\underline{t}),$
$HO(m,m;\underline{t})$ and $K(m,n;\underline{t})$ for any $m,n>2$,
$\underline{r}\in \mathbb{N}^{p+2}$ and $\underline{t}\in
\mathbb{N}^{m}.$
\end{proof}
\begin{remark}\label{remark5.4}
By  Corollaries \ref{c1} and  \ref{c2},
 the family of finite dimensional simple  special odd contact
superalgebras does contain ``strange" ones which are not isomorphic
to any simple Lie superalgebras of Cartan type $W, S, H, K, HO, KO$
or $SHO.$ A complete consideration for the isomorphism problem between modular Lie superalgebras of Cartan type is
 beyond the
scope of the present paper.
\end{remark}

\end{document}